\documentclass[12pt]{amsart}
\usepackage{latexsym}           %only needed for \leadsto in intro
\usepackage{amssymb}
\usepackage{amsxtra}
\usepackage{lscape}
\usepackage{epsfig}             % for \epsfig
\usepackage{amsfonts}
\newcommand\beq{\begin{equation}}
\newcommand\eeq{\end{equation}}

\newcommand{\IP}{\mathbb{P}}                                     
                           
\newcommand{\IR}{\mathbb{R}}                           
\newcommand{\IC}{\mathbb{C}}
\newcommand{\IZ}{\mathbb{Z}}

\newcommand{\M}{\mathcal{M}}

\newcommand{\cA}{\mathcal{A}}

\newcommand{\cC}{\mathcal{C}}
\newcommand{\bcC}{\boldsymbol{\mathcal{C}}}
\newcommand{\bcO}{\boldsymbol{\mathcal{O}}}

\newcommand{\F}{\mathcal{F}}

\newcommand{\cO}{\mathcal{O}}

\newcommand{\g}{       \mathfrak{g}     }

\newcommand{\lsl}{\mathfrak{sl}} 
\newcommand{\PVI}{$\text{P}_{\text{VI}}$}   %text for PVI
\newcommand{\bth}{{\bf \theta}}
\newcommand{\bm}{{\bf m}}
\newcommand{\bM}{{\bf M}}
\newcommand{\bi}{{\bf i}}
\newcommand{\bj}{{\bf j}}
\newcommand{\pf}{\begin{bpf}}

\newcommand{\pfms}{\begin{bpfms}}
\newcommand{\epf}{\end{bpf}\hfill$\square$\\}           % end proof
\newcommand{\epfms}{\end{bpfms}\hfill$\square$\\}               % end proof

\newcommand{\idea}{\begin{bidea}}

\newcommand{\eidea}{\end{bidea}\hfill$\square$\\}           % end proof

\newcommand{\sketch}{\begin{bsketch}}

\newcommand{\esketch}{\end{bsketch}\hfill$\square$\\}           % end proof

\newcommand{\wt}{\widetilde}

\newcommand{\al}{\alpha}
\newcommand{\alc}{\alpha\spcheck}
\newcommand{\be}{\beta}
\newcommand{\ga}{\gamma}
\newcommand{\de}{\delta}
\newcommand{\Ga}{\Gamma}
\newcommand{\si}{\sigma}
\newcommand{\om}{\omega}
\newcommand{\Si}{\Sigma}

\newcommand{\bk}{{\bf k}}

\newcommand{\Sym}{\text{\rm Sym}}

\newcommand{\re}{{\text{\rm Re}}}

\newcommand{\tr}{\text{\rm Tr}}
\newcommand{\Tr}{\text{\rm Tr}}
\newcommand{\Hom}{\text{\rm Hom}}
\newcommand{\Aut}{\text{\rm Aut}}

\newcommand{\SL}{\text{\rm SL}}
\newcommand{\PSL}{\text{\rm PSL}}
\newcommand{\GL}{\text{\rm GL}}
\newcommand{\PGL}{\text{\rm PGL}}

\newcommand{\SO}{\text{\rm SO}}

\newcommand{\SB}{\text{\rm SB}}

\newcommand{\Diff}{\text{\rm Diff}}
\newcommand{\Out}{\text{\rm Out}}

\newcommand{\mc}{\text{\rm M}_{0,4}}
\newcommand{\wmc}{\wt {\text{\rm M}}_{0,4}}

\newcommand{\Aff}{{\text{\rm Aff}}}

\newcommand{\sdp}{{\ltimes}}

\newcommand {\eps}{\varepsilon}

\newcommand{\spq}{/\!\!/}

\newenvironment{eqn}
        {\begin{equation}}
        {\end{equation}}

\theoremstyle{plain}
\newtheorem {hypo}{\bf\hspace{-\parindent}Hypothesis}

\newtheorem {thm}{Theorem}
\newtheorem {prop}[hypo]{Proposition}%[section]

\newtheorem {lem}[hypo]{Lemma}%[section]

\theoremstyle{definition}
\newtheorem {defn}[hypo]{Definition}%[section]
\theoremstyle{remark}
\newtheorem {rmk}[hypo]{Remark}%[section]

\begin{document}

%\date{\today}

%\title[Symplectic and Quasi-Hamiltonian Geometry]{Symplectic and
%Quasi-Hamiltonian Geometry of Meromorphic Connections over Riemann Surfaces}
%\title{From Klein to Painlev\'e via Fourier, Laplace and Jimbo}
\title[Icosahedral solutions]
{The fifty-two icosahedral solutions to Painlev\'e VI}
\author{Philip Boalch}
\address{\'Ecole Normale Sup\'erieure\\
45 rue d'Ulm\\
75005 Paris\\
France} 
\email{boalch@dma.ens.fr}


\begin{abstract}
The solutions of the (nonlinear) 
Painlev\'e VI differential 
equation having icosahedral linear monodromy group will
be classified up to equivalence under
Okamoto's affine $F_4$ Weyl group action 
and many properties of the solutions will be given.

There are $52$ classes,
the first ten of which correspond directly to
the ten icosahedral entries on Schwarz's list of algebraic solutions of the 
hypergeometric equation. 
The next nine solutions are simple deformations of known
\PVI\  
solutions (and have
less than five branches)
and five of the larger solutions  are already known, due to work of
Dubrovin and Mazzocco and Kitaev. 

Of the remaining $28$ solutions we 
will find $20$ explicitly using the method of \cite{k2p} 
(via Jimbo's asymptotic formula).
%and will show the others may be found, in principle, 
%by the same method
%(the largest solution would be  a degree $72$ polynomial representing a
%certain genus seven plane curve). 
Amongst those constructed there is one solution
that is ``generic'' in that its parameters lie on 
{\em none} of the affine $F_4$ hyperplanes, 
one that is equivalent to the Dubrovin--Mazzocco elliptic solution
and three elliptic solutions that are related to the Valentiner 
three-dimensional complex reflection group, the largest having $24$ branches.
\end{abstract}

%\dedicatory{}

%*************************************

\maketitle

\vspace{0.1cm}
{\footnotesize
\bf CONTENTS}
\begin{centering}
\footnotesize
$$
\begin{array}{lcr}
\ref{sn: intro}. \text{ Introduction} & \ \qquad\qquad\qquad\ &
\pageref{sn: intro} \\
\ref{sn: triples}. 
\text{ Generating triples} &&
\pageref{sn: triples} \\
\ref{sn: affine}. 
\text{ Parameter equivalence---affine Weyl groups} &&
\pageref{sn: affine} \\
\ref{sn: geom}. 
\text{ Geometric equivalence---mapping class groups} &&
\pageref{sn: geom} \\
\ref{sn: tables}. \text{ Properties of the $52$ solutions} &&
\pageref{sn: tables} \\
\ref{sn: generic}. \text{ The generic icosahedral solution} &&
\pageref{sn: generic} \\
\ref{sn: moreex}. \text{ More examples} &&
\pageref{sn: moreex} \\
\ref{sn: valentiner}. \text{ The Valentiner solutions} &&
\pageref{sn: valentiner} \\
\end{array}
$$
\end{centering}

%\nocite{FG}            %Springer LNM 1620 

\renewcommand{\baselinestretch}{1.02}            % 1.2 for distrib.
\normalsize
%\footnotesep=0.4cm             %0.4cm  with stretch 1.2 or 1.25 

\begin{section}{Introduction}\label{sn: intro}
\ 

The Painlev\'e VI equation ($\text{P}_{\text{VI}}$) 
is a second order nonlinear differential equation which governs 
the isomonodromic deformations of linear systems of differential equations of
the form
\beq\label{eq: lin syst}
\frac{d}{dz}-\left(
\frac{A_1}{z}+\frac{A_2}{z-t}+\frac{A_3}{z-1}\right),\qquad
A_i\in\g:=\lsl_2(\IC)
\eeq
as the second pole position $t$ varies in $\IP^1\setminus\{0,1,\infty\}$.
(The general case---varying all four pole positions---reduces to this case
using automorphisms of $\IP^1$.)

Painlev\'e VI is notoriously difficult to solve explicitly, and 
indeed it has been proved that `most' solutions are 
new transcendental functions.

Upon fixing the eigenvalues of the residues $A_i$ (and of the residue $A_4:=
-\sum_1^3A_i$ at infinity), and identifying two systems if they are related by
a constant gauge transformation, one obtains a moduli space
$$\bcO:= \cO_1\times\cdots \cO_4\spq G=
\left\{\ (A_1,\ldots A_4)\in \cO_1\times \cdots\cO_4 
\ \bigl\vert \ \text{$\sum A_i =0$}\ \right\}/G$$ 
of such systems, which is (complex) two-dimensional in general,
where $G=\SL_2(\IC)$ and $\cO_i\subset \lsl_2(\IC)$ is the adjoint orbit of
elements having the 
chosen eigenvalues, assumed non-resonant here.
(Traditionally one parameterises the choice of the four adjoint
orbits $\cO_i$ by four complex
numbers $\theta_i$ such that
$A_i \text{ has eigenvalues } \pm\theta_i/2.$)

Geometrically (see below) each $\text{P}_{\text{VI}}$ equation 
amounts to a (nonlinear) connection on the
trivial fibre bundle
$$\M^*:=\bcO\times B \to B$$
where the base $B:=\IP^1\setminus\{0,1,\infty\}$ is the domain of $t$.

Thus, roughly speaking, the set 
\begin{equation}\label{eq: big poisson var}
\g^4\spq G = \left\{\ (A_1,\ldots A_4)\in \g^4 
\ \bigl\vert \ \sum A_i =0\ \right\}/G
\end{equation}
of residues is foliated by a family of surfaces $\bcO$ parameterised by 
$\theta_1,\ldots\theta_4$ (i.e. as the $\theta$'s vary 
the surfaces $\bcO\subset \g^4\spq G $ sweep out most of $\g^4\spq G$).

Now in \cite{OkaPVI} K. Okamoto has
defined, again roughly speaking, a birational action of the affine
Weyl group $W_a(F_4)$ on the space of systems
$$\g^4\spq G\times B.$$
%of the space \eqref{eq: big poisson var} and the base $B$.
(A point of this product corresponds to a choice of residues $A_i$ and a choice
of $t\in B$, and so determines a system as in \eqref{eq: lin syst}.)
The action on $B$ is via M\"obius transformations $\mu$ permuting $0,1,\infty$
and, as we will confirm, the action on  $(\theta_1,\ldots,\theta_4)\in \IC^4$
is the standard $W_a(F_4)$ action.
Two key properties of Okamoto's action are:

$\bullet$
It maps each leaf $\bcO\times B\subset \g^4\spq G\times B$ to another leaf,
say $\bcO'\times B$, and

$\bullet$
It relates the $\text{P}_{\text{VI}}$ connection on  $\bcO\times B$ to the 
$\text{P}_{\text{VI}}$ connection on 
$\bcO'\times B$ (so that local $\text{P}_{\text{VI}}$ solutions $s:U\to \bcO$ 
for $U\subset B$
map to local solutions $s':\mu(U)\to \bcO'$).

The author's prior understanding of the $\text{P}_{\text{VI}}$ folklore
 was that the classical
solutions of $\text{P}_{\text{VI}}$ (i.e. those which are not `new' 
transcendental functions) had
$\theta$-parameters lying on one or more of the reflection hyperplanes of the
$W_a(F_4)$ action (the idea being that such solutions had more symmetry).

One of the main aims of this paper is to show that this is not necessarily the
case---an explicit algebraic solution will be written down with parameters on
{\em none} of the affine $F_4$ hyperplanes.

To explain the strategy used let us first recall more of the geometrical path
to $\text{P}_{\text{VI}}$.

The monodromy representation of a linear system \eqref{eq: lin syst}
is a point of the space
$$\bcC=\Hom_\cC(\pi_1(\IP^1\setminus\{0,t,1,\infty\}),G)/G$$
of conjugacy classes of  representations of the 
fundamental group of the four-punctured sphere, where the representations are
restricted to take the simple loop around the $i$th puncture
into the conjugacy class $\cC_i:=\exp(2\pi\sqrt{-1} \cO_i)\subset G$.
Upon choosing appropriate loops generating 
$\pi_1(\IP^1\setminus\{0,t,1,\infty\})$, 
$\bcC$ is identified with
\beq\label{eq: explicit C}
\left\{ (M_1,\ldots,M_4)\ \bigl\vert\ 
M_i\in \cC_i,\  M_4\cdots M_1=1\ \right\}/G
\eeq
the multiplicative analogue of $\bcO$, and is similarly seen to be
two-dimensional in general.

As the position $t$ of the second pole varies these surfaces $\bcC$ fit
together into a fibre bundle 
$$M \to B,$$
the fibre over $t\in B$ being the surface $\bcC$ associated to the
four-punctured sphere $\IP^1\setminus\{0,t,1,\infty\}$.

This fibre bundle $M\to B$ has a natural complete flat connection  
on it (in other words it is a local system of varieties):
as $t\in B$ is varied slightly we can identify two 
nearby fibres of $M$ by using
the same loops to identify both fibres with \eqref{eq: explicit C}
and therefore with each other.

The $\text{P}_{\text{VI}}$ equation is obtained by pulling back this connection on $M\to B$ along
the relative Riemann--Hilbert map 
$$\nu:\M^* \to M$$
(taking systems to their monodromy).
Choosing particular ($t$-dependent) coordinates on the fibres of $\M^*$,
writing out what one gets and eliminating one coordinate, yields $\text{P}_{\text{VI}}$:
$$\frac{d^2y}{dt^2}=
\frac{1}{2}\left(\frac{1}{y}+\frac{1}{y-1}+\frac{1}{y-t}\right)
\left(\frac{dy}{dt}\right)^2
-\left(\frac{1}{t}+\frac{1}{t-1}+\frac{1}{y-t}\right)\frac{dy}{dt}  $$
$$
%\quad\qquad\qquad
\quad\ +\frac{y(y-1)(y-t)}{t^2(t-1)^2}\left(
\al+\be\frac{t}{y^2} + \gamma\frac{(t-1)}{(y-1)^2}+
\delta\frac{ t(t-1)}{(y-t)^2}\right) $$
where the constants $\al,\be,\ga,\de$ 
are related to the $\theta$-parameters as follows:
\begin{equation}\label{thal}
\al=(\theta_4-1)^2/2, \qquad \be=-\theta_1^2/2, \qquad 
\ga=\theta^2_3/2, \qquad \delta=(1-\theta_2^2)/2.
\end{equation}

From this picture it is immediate that the branching of solutions to $\text{P}_{\text{VI}}$
(horizontal sections of the connection on $\M^*$)
corresponds to the monodromy of the connection on $M$.
But the connection on $M$ is complete and so its monodromy amounts to an
action of the fundamental group $\F_2=\pi_1(B)$ of the base $B$ on a fibre
$\bcC$. This action of $\F_2$ (the free group on two generators) 
can be written down explicitly in terms of the standard Hurwitz
braid group action.

The simplest solutions of $\text{P}_{\text{VI}}$ should be those with only a finite number of
branches. For example if we take a linear 
system \eqref{eq: lin syst} which has a
basis of algebraic solutions then the corresponding $\text{P}_{\text{VI}}$ solution (controlling
the isomonodromic deformations of \eqref{eq: lin syst})
will be finite branching.
One way to see this is to recall that \eqref{eq: lin syst} will have a basis
of algebraic solutions if and only if its monodromy group
$\langle M_1,M_2,M_3 \rangle \subset G$ is finite, i.e. the $M_i$ are
generators of a finite subgroup $\Ga$ of $G=\SL_2(\IC)$.
Now the triple $M_1,M_2,M_3$ represents a point of a surface $\bcC$ and the
$\F_2$ action on $\bcC$ (the monodromy of the $\text{P}_{\text{VI}}$ solution) acts within the set
of triples of generators of $\Ga$.
Thus the number of branches of the $\text{P}_{\text{VI}}$ solution is bounded, 
e.g. by $\bigl\vert \Ga\bigl\vert^3$.
The idea of looking for solutions of $\text{P}_{\text{VI}}$ starting from a finite subgroup of
$\SL_2(\IC)$
goes back at least to Hitchin \cite{Hit-Poncelet} (see also \cite{Hit-Octa}).

In this paper we will examine the set of $\text{P}_{\text{VI}}$ solutions which arise in this way
upon taking $\Ga$ to be the binary icosahedral group.

%\begin{thm}
\vspace{0.2cm}
\noindent{\bf Theorem A.} (Icosahedral Classification)\ {\em
Upto equivalence under Okamoto's $W_a(F_4)$ action, 
there are precisely $52$ solutions to $\text{P}_{\text{VI}}$ 
having icosahedral linear monodromy
group. The possible genera are $0,1,2,3,7$ and the largest solution has $72$
branches.    }
\vspace{0.2cm}
%\end{thm}

%We find there are $52$ equivalence classes of icosahedral solutions
%with possible genera $0,1,2,3,7$. 

More details about the solutions can be found in Table $1$.
Examining the parameters of the solutions we find
exactly one such solution turns out to have generic parameters,
corresponding to 
choosing one of the four generators $M_1,M_2,M_3,M_4$
to be in (a lift to $\SL_2(\IC)$ of) 
each of the four nontrivial conjugacy classes of $A_5$, the icosahedral
rotation group. (The reader may like to confirm that 
there are such four-tuples of
elements of $A_5$ having product the identity.)  This leads to:

%\begin{thm}\label{thm: exists genc}
\vspace{0.2cm}
\noindent{\bf Theorem B.} (Generic Solution)\ {\em
There is an algebraic solution to the sixth Painlev\'e equation whose
parameters lie on none of the reflecting hyperplanes of Okamoto's affine $F_4$
(or $D_4$) action.}
\vspace{0.2cm}
%\end{thm}

\noindent
{\bf Proof.}\quad
Set 
$(\theta_1,\theta_2,\theta_3,\theta_4)=(2/5, 1/2, 1/3, 4/5)$ 
and consider the rational functions 
$$
y=-\frac{9 s (s^2+1) (3 s-4) (15 s^4-5 s^3+3 s^2-3 s+2)}
{(2 s-1)^2(9 s^2+4)(9 s^2+3 s+10)},\quad 
t=\frac{27 s^5 (s^2+1)^2 (3 s-4)^3}{4(2 s-1)^3(9 s^2+4)^2}\quad
\begin{matrix} \\ \\ \square \end{matrix}
$$

\ 

The main technical tool used in the construction of this solution is the
precise 
formula of M. Jimbo (see \cite{Jimbo82} and \cite{k2p} Theorem 4) for the 
leading term in the asymptotic expansion at zero of generic
$\text{P}_{\text{VI}}$ solutions. 
In brief, using the $\text{P}_{\text{VI}}$ equation 
these leading terms determine the Puiseux expansions of 
each branch of the
solution at zero and, taking sufficiently many terms, 
these determine the solution completely since it is algebraic.

One of the basic facts the author came to appreciate during 
the construction
of this generic solution 
is that even though two solutions may be equivalent by Okamoto
transformations, the size of the polynomial defining 
them may well vary dramatically;
one should try to choose the parameters for which they become as simple as
possible (which is still something of an art). 
For example the first solution curve $F(y,t)=0$ found, defining an
equivalent generic solution, took a page to write down and involved about 
a hundred 
 twenty-digit integers. A more perspicacious choice of parameters reduced
the size of the polynomial $F$ and enabled the above parameterisation to be
computed. 

This led to the question of whether 
there is a better choice of equivalent parameters for 
the elliptic icosahedral solution of Dubrovin and Mazzocco \cite{DubMaz00},
for which the solution curve $F(y,t)$ took about ten pages 
to write down (in the
preprint version of op. cit. on the math archive), and for which a
parameterisation was not possible to compute.

%\begin{thm}\label{thm: dub maz equiv}
\vspace{0.2cm}
\noindent{\bf Theorem C.}\ {\em
The Dubrovin--Mazzocco elliptic solution is equivalent to the
solution with parameters 
%$(\theta_1,\theta_2,\theta_3,\theta_4)=(1/3, 1/3, 1/3, 1/3)$ 
$\theta_i=1/3$, $i=1,2,3,4$
given by the
functions $y,t$ on the elliptic curve $$u^2=s\,(8\,s^2-11\,s+8)$$
\noindent
where 
$$
y=
\frac{1}{2}-
{\frac {8\,{s}^{7}-28\,{s}^{6}+75\,{s}^{5}+31\,{s}^{4}-269\,
{s}^{3}+318\,{s}^{2}-166\,s+56}
{18\,u \left( s-1 \right)  \left( 3\,{s}^{3}-4\,{s}^{2}+4\,s+2 \right) }}$$

$$t=
\frac{1}{2}+
{\frac { \left( s+1 \right)  
\left( 
32\,({s}^{8}+1)-320\,({s}^{7}+s)+1112\,({s}^{6}+s^2)
-2420\,({s}^{5}+s^3)+3167\,{s}^{4}
\right) }
{54\,{u}^{3}\,s\, \left( s-1 \right) }}.
$$
In particular this elliptic curve is birational to that defined by the
$10$-page polynomial.
}
\vspace{0.2cm}
%\end{thm}

We remark that this solution was constructed directly with the above parameter
values, rather than by transforming the curve of Dubrovin and Mazzocco.
(It is now straightforward to apply Okamoto transformations to the above 
parameterised solution to obtain a parameterisation of their curve; $t,u,s$
remain unchanged but $y$ becomes more complicated: see 
Section \ref{sn: tables}.)

The author was first motivated to examine the set of icosahedral solutions 
to \PVI\  for the following
reason.
In previous articles \cite{pecr,k2p} the author studied 
an alternative realisation where \PVI\  governs the isomonodromic
deformations of certain rank three Fuchsian systems having four poles 
on $\IP^1$ 
and having rank one residues at three of the poles.
In particular 
\cite{k2p} explained how to relate these systems, and their monodromy data, 
to the standard $\SL_2$ viewpoint described above.
In terms of monodromy representations this yields a direct way 
to construct finite $\F_2$ 
orbits of $\SL_2(\IC)$ triples, from each triple of reflections generating a
(finite) three-dimensional complex reflection group. 
For example the Klein solution was constructed in \cite{k2p}
starting from a triple of generating reflections
of the smallest non-real exceptional three-dimensional
complex reflection group, the Klein group.
In that case the corresponding $\SL_2(\IC)$ triple is not equivalent
to that of any finite subgroup of $\SL_2(\IC)$.

Somewhat disappointingly most of the other three dimensional complex
reflection groups seem to lead to known Painlev\'e VI solutions (or simple
deformations of them).
However the largest exceptional three-dimensional
complex reflection group, the Valentiner group (a six fold cover of $A_6$ of
order $2160$), does yield new solutions: one finds there are three
inequivalent triples of generating reflections, with $\F_2$ orbits of sizes 
$15,15,24$ respectively, all corresponding to genus $1$ solutions. 
However the corresponding triples in $\SL_2(\IC)$ all turn out to generate the
binary icosahedral group. (In particular this gives an unexpected
relationship between $A_6$ and $A_5$.)
Thus we realised there are other interesting icosahedral solutions distinct
from those previously found, and so became curious to see any others
that might occur.
These Valentiner solutions now appear as 
rows $37, 38, 46$ of Table $1$, and we have managed to construct all three 
solutions explicitly. 
(Currently the $24$ branch solution is the highest degree explicit 
algebraic solution to \PVI.)

\vspace{0.2cm}
\noindent{\bf Theorem D.} (Valentiner Solutions)\ {\em
There are three inequivalent triples of reflections generating the Valentiner
complex reflection group having $\F_2$ orbits of sizes $15,15,24$
respectively. 
The corresponding \PVI\ solutions all have genus one and are equivalent to
icosahedral solutions. (They will appear in Section \ref{sn: valentiner}.)} 
\vspace{0.2cm}

Following the procedure of \cite{k2p} these solutions give explicit families
of rank three, four-poled Fuchsian systems
having monodromy
the Valentiner group, generated by reflections.

\ 

The layout of this article is as follows.
Section \ref{sn: triples} 
describes convenient parameters on, and enumerates, the set $S$ of
conjugacy classes of triples of generators of $\Ga$. 
Section \ref{sn: affine} then counts the set of $W_a(F_4)$ orbits of the
$\theta$-parameters that arise from $S$.
Since $\Ga$ is finite all such $\theta$'s are real 
 and so this amounts to reflecting the parameters 
into a fundamental domain (the closure of
 an alcove) for the real action of $W_a(F_4)$.
This gives a lower bound, $52$, on the number of inequivalent icosahedral
solutions.  
Section \ref{sn: geom} then proves $52$ is also an upper bound by examining
the natural action of the mapping class group (and of the centre of $G^3$) on
triples of generators of $\Ga$, and relating this to Okamoto's $W_a(F_4)$
action. Combined with Section \ref{sn: affine} this gives the desired
classification. 
Section \ref{sn: tables} then lists many properties of the $52$ solutions 
(and describes the relation 
with Schwarz's list, with Lam\'e equations
and with 
the previously published icosahedral solutions of Dubrovin and Mazzocco and
Kitaev). 
Section \ref{sn: generic} discusses the generic icosahedral solution,
then section
\ref{sn: moreex}
presents some other explicit 
icosahedral solutions, including all the outstanding
genus zero solutions as well as some others of genus one.
(Remark \ref{rmk: qt}, added after the rest of this paper was written, explains
how the remaining solutions may be obtained.)
Finally 
section \ref{sn: valentiner} presents the Valentiner solutions, which was our
starting point.\footnote{
To aid the reader who is interested in examining the solutions of this paper
(for example to draw the corresponding dessins d'enfants),
a Maple text file of the solutions has been included with the source file on
the math arxiv (math.AG/0406281). 
This may be downloaded by clicking on ``Other formats'' and
unpacked with the commands `gunzip 0406281.tar' and `tar -xvf 0406281.tar', 
at least on a Unix system.
}

Similar considerations may also be applied to the cases of the 
tetrahedral and octahedral
groups; this has now been done (and all such \PVI\ solutions are now known
explicitly). Details will appear elsewhere \cite{octa}.

\

\Small
{\bf Acknowledgments.}\ \ 
The author is grateful to  
Patrick X. Gallagher for pointing out Hall's paper \cite{Hall36},
to Anton Alekseev for inviting him to visit
the University of Geneva Mathematics  Department
where the first version of the generic solution was found
(supported by the Swiss NSF), and especially to
Mark van Hoeij (both for direct help finding the more difficult rational
parameterisations and indirectly for having written the Maple algebraic curves
package).
Most of the computations for this article were done on the 
Medicis computers at
\'Ecole Polytechnique. 
The Valentiner solutions were found whilst the
author was visiting Kobe University Mathematics Department, December 2004.

\end{section}

\begin{section}{Generating Triples}\label{sn: triples}

Let $G=\SL_2(\IC)$ and consider the binary icosahedral group 
$\Ga\subset G$ of order $120$.
It has a center of order $2$ and the quotient 
$\Ga/\pm \subset \PSL_2(\IC)\cong \SO_3(\IC)$ is the icosahedral
group $A_5$.
In terms of unit quaternions, explicit generators of $\Ga$ are 
\cite{conway-smith}:
$$
\Ga = \left\langle \ (-1+\bi+\bj+\bk)/2,\ (\bi+\si\bj+\tau\bk)/2\ 
\right\rangle
$$ 
%\frac{1}{2}
where $\si=(\sqrt{5}-1)/2, \tau=(\sqrt{5}+1)/2$
and $
\bi=\left(\begin{smallmatrix}
i & 0  \\
0 & -i
\end{smallmatrix}\right),
\bj=\left(\begin{smallmatrix}
0 & 1  \\
-1 & 0
\end{smallmatrix}\right),
\bk=\left(\begin{smallmatrix}
0 & i  \\
i & 0
\end{smallmatrix}\right).
$

Our first aim is to study the set of triples of generators of $\Ga$.
Suppose we have a triple $M_1,M_2,M_3\in \Ga$ which generate $\Ga$ (rather than
a proper subgroup). Denote this triple by $\bM$:
$$\bM= (M_1,M_2,M_3).$$
Then define $M_4\in \Ga$ by the requirement that 
\begin{equation} \label{eq: monod. relation}
M_4M_3M_2M_1 = 1
\end{equation}
and consider the seven-tuple of numbers
$$\bm=\bm(\bM)=(m_1,m_2,m_3,m_4,m_{12},m_{23},m_{13})$$
where
$$m_i:=\tr(M_i),\qquad m_{ij}:=\tr(M_iM_j).$$

\begin{lem}
Two triples of generators of $\Ga$ are conjugate (in $G$ or $\Ga$) 
if and only if they have the same
seven-tuple $\bm$.
\end{lem}
\sketch
This follows since these traces generate the ring of invariants of the
diagonal conjugation action of $G$ on $G^3$, and that $\Ga$ is its own
normaliser in $G$. 
\esketch

There is in fact a formula to count the conjugacy classes of generating triples:
\begin{lem}[\cite{Hall36}]
There are $26688$ conjugacy classes of triples of generators of $\Ga$.
\end{lem} 
\pf
P. Hall shows (\cite{Hall36} p.146) that the number of $n$-tuples
of generators of $\Ga$ is
$$120^n - 5 (24^n)-6 (20^n)-10 (12^n) +20 (6^n) +60 (4^n) - 60 (2^n).$$
%120^n - 5 *(24^n)-6* (20^n)-10* (12^n) +20 *(6^n) +60 *(4^n) - 60 *(2^n);
Setting $n=3$ and dividing by $60$ (since the center of $\Ga$ acts trivially
when conjugating) gives the result.
\epf

Thus the set
$$S:= \{ \bm(\bM)\ \bigl\vert\  \langle M_1,M_2,M_3\rangle=\Ga \ \}$$
 of seven-tuples of invariants of generating triples of $\Ga$, has
cardinality $26688$.
Fortunately there are quite strong notions of equivalence for elements of $S$, 
which will dramatically reduce this number.
In the next sections 
we will define and study two equivalence relations (parameter equivalence
 and geometric equivalence) on $S$,
which will turn out to have the same $52$ equivalence classes.
\end{section}

\begin{section}{Parameter equivalence---affine Weyl groups}\label{sn: affine}

Given a seven-tuple $\bm\in S$ we can associate four parameters
$$ \bth= \bth(\bm):= (\theta_1,\theta_2,\theta_3,\theta_4)\in\IR^4$$
where $\theta_j\in\IR$ is determined from $\bm$ via:
$$m_j=2 \cos(\pi \theta_j),\qquad 0\leqslant\theta_j \leqslant 1$$
so that the corresponding 
matrix $M_j\in \Ga$ has eigenvalues $\{\exp(\pm \pi i \theta_j)\}$.

\begin{defn}
Two seven-tuples $\bm,\bm'$ are {\em parameter equivalent} if
their parameters $\bth,\bth'$ are in the same orbit of the standard action of
the affine Weyl group of type $F_4$ on $\IR^4$.
\end{defn}

In order to explain this,
let us 
briefly recall some basic facts about root
systems and the corresponding affine Weyl group actions
(for more details see e.g. Bourbaki \cite{bbkLie}). 

Let $V$ be a real four-dimensional Euclidean vector space with orthonormal
basis $\eps_1,\ldots,\eps_4$.
The Euclidean inner product will be denoted $(u,v)$ and used to identify $V$
with its dual $V^*$. 
 Let $O(V)$ denote the group of linear transformations of $V$ preserving the
inner product and let $\Aff(V)\cong O(V)\sdp V$ 
denote the group of affine Euclidean transformations of $V$
(i.e. those of the form $v\mapsto g(v)+w$ for some $g\in O(V), w\in V$). 
A vector in $V$ will be denoted $\sum \theta_i\eps_i$ with
$\theta_i\in \IR$ (the indices on $\eps_i$ and $\theta_i$ will always run from
$1$ up to $4$).

The standard $F_4$ root system is the following set of $48$ vectors in $V$:
$$F_4=\{ \ \pm \eps_i, \ \pm \eps_i\pm \eps_j (i<j),\  
(\pm\eps_1\pm\eps_2\pm\eps_3\pm\eps_4)/2 \}.$$

Each root $\alpha\in F_4$ determines a coroot 
$\alc = \frac{2\al}{(\al,\al)}$ as well as
a hyperplane $L_\al$ in $V$:
$$L_\alpha:= \{\  v\in V\ \bigl\vert\ (\alpha,v)=0\ \}.$$
In turn $\al$ 
determines an orthogonal reflection $s_\alpha$, the reflection in this
hyperplane:
$$s_\alpha(v)= v- 2\frac{(\al,v)}{(\al,\al)}\al=v- (\alc,v)\al.$$

The Weyl group $W(F_4)\subset O(V)$ is the group 
generated by these reflections:
$$W(F_4)= \langle\ s_\al\ \bigl\vert\ \al\in F_4\ \rangle$$
which is of order $1152$.

Similarly the choice of a root $\alpha\in F_4$ and an integer $k\in \IZ$
determines an affine hyperplane $L_{\al,k}$ in $V$:
$$L_{\alpha,k}:= \{\  v\in V\ \bigl\vert\ (\alpha,v)=k\ \}$$
and the reflection $s_{\alpha,k}$ 
in this hyperplane is an affine Euclidean transformation 
$$s_{\al,k}(v)= s_\al(v)+k\alc.$$

The affine Weyl group $W_a(F_4)\subset \Aff(V)$ is the group 
generated by these reflections:
$$W_a(F_4)= \langle\ s_{\al,k}\ \bigl\vert\ \al\in F_4,\  k\in\IZ\ \rangle$$
which is an infinite group isomorphic to the semi-direct product
of $W(F_4)$ and the coroot lattice $Q(F_4\spcheck)$ (which is the lattice in
$V$ generated by the coroots $\alc$).

Now in \cite{OkaPVI} Section 3, Okamoto defines a birational
action of (a copy of) $W_a(F_4)$ on a $7$-dimensional 
space of linear differential equations,
the ``total phase space'' of 
Painlev\'e VI,
(involving the four parameters plus the canonical coordinates $p,q$ and the
time variable $t$---this is essentially the space $\g^4\spq G\times B$ 
in the introduction, with $p,q$ being coordinates on $\bcO$;
given an isomonodromic family of linear equations it is the function $y=q$
which solves $\text{P}_{\text{VI}}$).
This action descends to an action on just the space of the four parameters 
(denoted $v_i$ in \cite{OkaPVI}).
By relating Okamoto's four parameters to the $\theta$-parameters used here
we see that Okamoto defines an embedding 
$\iota:W_a(F_4)\hookrightarrow\Aff(V)$.
\begin{lem}
Okamoto's embedding maps his copy of $W_a(F_4)$ isomorphically onto the
standard $W_a(F_4)\subset \Aff(V)$.
\end{lem}
\pf
The action of $\iota(W_a(F_4))$ is generated (\cite{OkaPVI} p.364) by the
reflections in the five 
hyperplanes bounding the alcove:
\begin{equation}\label{eq: ok alc}
v_2> v_3 > v_4 > 0,\qquad v_1> v_2+v_3+v_4,\qquad v_1+v_2< 1
\end{equation}
where $v_1=\theta_3-1, v_2=\theta_1, v_3=\theta_2, v_4=\theta_4-1$,
whereas the standard $W_a(F_4)$ is generated by the reflections in the
hyperplanes bounding the standard alcove:
\begin{equation}\label{eq: std alc}
\theta_2> \theta_3 > \theta_4 > 0,\qquad 
\theta_1> \theta_2+\theta_3+\theta_4,\qquad 
\theta_1+\theta_2< 1.
\end{equation}
One may show $\iota(W_a(F_4))\subset W_a(F_4)$ by finding $g\in W_a(F_4)$
mapping \eqref{eq: ok alc} isomorphically onto \eqref{eq: std alc}.
(Such $g$ may be found by applying the procedure of Proposition   
\ref{prop: at least 52} below to an interior point of \eqref{eq: ok alc}.)
Similarly for the reverse inclusion.
\epf 

\begin{rmk}
The reason we are being careful here and speaking of different copies of 
$W_a(F_4)$, is that the analogous result is not true for Okamoto's affine $D_4$
action. 
Recall from \cite{OkaPVI} that Okamoto starts by defining an action of 
$W_a(D_4)$ (which fixes the time variable $t$) and the action of $W_a(F_4)$ 
is obtained by adding some more generators. However when written as an action
on $V$ (our space of $\theta$'s) $W_a(D_4)$ is not embedded in $\Aff(V)$ as
the standard $W_a(D_4)$, but rather as $W_a(D^-_4)$ where
$D_4^-$ is the set of $24$ {\em short} roots of $F_4$:
$$D_4^-:= \{\ \pm \eps_i,\  (\pm\eps_1\pm\eps_2\pm\eps_3\pm\eps_4)/2 \ \},$$
whereas the standard $D_4$ is the set of long roots of $F_4$.
Moreover one then naturally has that $W_a(F_4)$ is the normaliser of 
$W_a(D^-_4)$ in $\Aff(V)$, and is an extension by $S_4$, the automorphisms of
the extended $D_4$ Dynkin diagram;
each $W_a(D^-_4)$ alcove is partitioned into 
$24=|S_4|$ copies of the $W_a(F_4)$ alcove. 
\end{rmk}

In this paper we are viewing two solutions to $\text{P}_{\text{VI}}$ as equivalent if they are
related by Okamoto's action of $W_a(F_4)$. 
Thus it is immediate that the $2\times 2$ linear monodromy data of 
any two equivalent $\text{P}_{\text{VI}}$ solutions will be parameter equivalent.
Hence by computing the set of $W_a(F_4)$ orbits of the 
set $\bth(S)\subset V$ we obtain a lower bound for the set of inequivalent 
icosahedral $\text{P}_{\text{VI}}$ solutions:

\begin{prop}\label{prop: at least 52}
There are at least $52$ inequivalent icosahedral solutions to $\text{P}_{\text{VI}}$.
\end{prop}
\pf
Direct computation---the 
standard procedure for computing affine Weyl group orbits is as follows.
The set of affine hyperplanes partitions $V$ into a set of disconnected
pieces, the alcoves,
and the affine Weyl group acts simply transitively on the set of these alcoves.
Every affine Weyl group orbit intersects the closure of any alcove 
in exactly one point.

Thus we choose an alcove $\cA$, so that $V/W_a(F_4)\cong \overline\cA$, 
and for each
point of $\bth(S)$ we find the corresponding point of $\overline\cA$. 
This is done by
repeatedly reflecting in the hyperplanes bounding $\cA$ until all the
inequalities determining $\overline\cA$ 
are satisfied. (This procedure will always 
terminate after a finite number of steps.)
Proceeding in this way we find (using Maple) that the set $\bth(S)$ leads to
precisely $52$ points of $\overline\cA$.
\epf

\begin{rmk}
%The author has yet to decide on his preferred choice of alcove; 
There  is clearly lots of choice of alcove---one would
perhaps eventually like
to find an alcove in which the corresponding $\text{P}_{\text{VI}}$ solutions have as simple
form as possible.
%, or to find an alcove in which has the largest possible
%intersection with $\bth(S)$. 
In the meantime we can use for example the standard alcove
\eqref{eq: std alc} or the alcove \eqref{eq: ok alc} suggested by 
Okamoto's work 
or that suggested by Noumi--Yamada's article \cite{NY-so8}:
$$
\al_2 > 0,\qquad \al_0>\al_1>\al_4>\al_3>0 
$$
where 
$\al_0=\theta_2, \al_1=\theta_4-1, \al_3=\theta_3, \al_4=\theta_1, 
\al_2=(1-\al_0-\al_1-\al_3-\al_4)/2$
(which is often convenient because the full birational 
action of $W_a(D_4)$ is given succinctly in \cite{NY-so8} and  
its extension to $W_a(F_4)$ is written in these terms in 
\cite{Masuda-PVI} 7.14).
\end{rmk}

Next we will look for a sharp 
upper bound on the number of icosahedral solutions.

\end{section}

\begin{section}{Geometric equivalence---mapping class groups}\label{sn: geom}

Let $a_1,a_2,a_3,a_4$ be four distinct points of the real 
two-dimensional sphere $S^2$ 
(say $a_1,\ldots a_4=0,\frac{1}{2},1,\infty$) and
consider the mapping class group of the sphere preserving the set of these
points:
$$\mc:= \pi_0(\Diff(S^2,\{a_1,a_2,a_3,a_4\}))$$
which is the {\em group} of connected components of the group of orientation 
preserving diffeomorphisms $f:S^2\to S^2$ such that  
$f(\{a_1,a_2,a_3,a_4\})=\{a_1,a_2,a_3,a_4\}$.

The following facts about $\mc$ will be useful:

$\mc$ is generated by elements $\om_i,\  i=1,2,3$ where $\om_i$ is related to
the Dehn twist swapping $a_i,a_{i+1}$ in an anti-clockwise sense. They 
satisfy the relations
given in \cite{Birman} p.164.

By mapping $\om_i$ to the permutation $(i,i+1)\in S_4$, one
obtains
the exact sequence
\begin{equation}\label{eq: mc to s4}
1\to \F_2\to \mc \to S_4 \to 1
\end{equation}
where the kernel (the pure mapping class group) is isomorphic to the free
group $\F_2$ on two letters, freely generated by $\om_1^2, \om_2^2$. 

Write $S^*=S^2\setminus\{a_1,a_2,a_3,a_4\}$ for the four-punctured sphere.
There is a natural map $\mc\to\Out(\pi_1(S^*))$ 
to the group of outer automorphisms (the group of all automorphisms modulo the
inner automorphisms) of the fundamental group of the
four-punctured sphere.
This is defined as follows: Given $f\in \mc$ one obtains an isomorphism 
$f_*:\pi_1(S^*)\to\pi_1(S^*)$, however the basepoint may well move, so one
needs to quotient by inner automorphisms.

This map induces an action of $\mc$ on the set of conjugacy classes of
representations of the
fundamental group of the four-punctured sphere:
$$\mc\to \Aut(\Hom(\pi_1(S^*),G)/G).$$
Explicitly the generators act as follows.
First choose simple positive loops $\ga_i$  around $a_i$ generating 
$\pi_1(S^*)$ such that $\ga_4\circ\cdots\circ \ga_1$ is contractible,
and let $M_i=\rho(\ga_i)\in G$ for any representation 
$\rho\in\Hom(\pi_1(S^*),G)$.
Then $\om_i$ fixes $M_j$ for $j\ne i,i+1, \ (1\leqslant j\leqslant 4)$ and
$$\om_i(M_i,\ M_{i+1})= (M_{i+1},\ M_{i+1}M_{i}M^{-1}_{i+1}).$$ 
In terms of the traces $m_i=\Tr(M_i),\ m_{ij}=\Tr(M_iM_j)$
generating the ring of $G$-invariant functions
on $\Hom(\pi_1(S^*),G)$ one finds (as in \cite{k2p} Lemma 1): 
\begin{gather}
\om_1(\bm)=
(m_2,\ m_1,\ m_3,\ m_4,\ m_{12},\ m_2 m_4 + m_1 m_3-m_{13}-m_{12}m_{23},\
m_{23}) \notag \\
\om_2(\bm)=
(m_1,\ m_3,\ m_2,\ m_4,\ m_{13},\  m_{23},\ m_3 m_4 + m_1
m_2-m_{12}-m_{23}m_{13}) \notag \\
\om_3(\bm)=
(m_1,\ m_2,\ m_4,\ m_3,\ m_{12},\  m_2 m_4 + m_1 m_3-m_{13}-m_{12}m_{23},\
m_{23})  \notag
\end{gather}
where $\bm=(m_1,\ m_2,\ m_3,\ m_4,\ m_{12},\ m_{23},\ m_{13})$.
(In computing this action 
we follow the conventions of \cite{Birman}; in the conventions
used in \cite{k2p} one has $\om_1=\be_2^{-1},\om_2=\be_1^{-1}$.) 

Let $t_0=1/2$ be a basepoint of $B=\IP^1\setminus\{0,1,\infty\}$ 
and choose loops $w_1$ (resp. $w_2$) based at $t_0$ encircling
$0$ (resp. $1$) once in a positive sense.
These two loops generate $\pi_1(B)\cong \F_2$ and there is a canonical map
$\pi: \pi_1(B) \to \mc$ mapping $\pi_1(B)$ isomorphically onto the pure
mapping class group, $\pi(w_i)=\omega_i^2$ for  $i=1,2$.
The action of this $\F_2$ corresponds to the nonlinear monodromy of 
Painlev\'e VI.
(To define $\pi$ geometrically, first define a map 
$B\to ((\IP^1)^4\setminus\text{diagonals})/S_4$
by mapping $t$ to the unordered set $\{0,t,1,\infty\}$.
Taking fundamental groups gives a map $\pi_1(B)\to \SB_4$  to the four-string
spherical braid group \cite{Birman} p.34.
Then recall $\mc$ is naturally the quotient of $\SB_4$ by its centre \cite{Birman}
Theorem 4.5.
The relation between the generators is described \cite{Birman} p.165.)

As well as the mapping class group there is another symmetry group acting on
the monodromy data that we wish to consider. 
Recall that there are precisely two connected Lie groups with Lie algebra
$\lsl_2(\IC)$: the simply connected group $\SL_2(\IC)$ with center $\pm 1$
and its quotient $\SL_2(\IC)/\pm 1 =\PSL_2(\IC)\cong \SO_3(\IC)$.
Thus any triple $\bM=(M_1,M_2,M_3)$  projects to a triple of elements of
$\PSL_2(\IC)$. We will say two triples $\bM, \bM'$ are {\em sign equivalent}
if they project to the same triple in $\PSL_2(\IC)$.
Said differently, let
$$\Si = \{ \ (\epsilon_1,\epsilon_2,\epsilon_3,\epsilon_4)\ \bigl\vert\ 
\epsilon_i=\pm 1,\ \Pi \epsilon_i=1\ \} \cong (\IZ/2)^3$$
be the group of even four-tuples of signs
then, since $M_4M_3M_2M_1=1$, we are acting on the four-tuple
$M_1,M_2,M_3,M_4$ with $\Si$ in the obvious way $M_i\mapsto \epsilon_iM_i$.

The mapping class group $\mc$ acts on $\Si$ via the map \eqref{eq: mc to s4}
to $S_4$ and the obvious action of $S_4$ permuting the $\epsilon_i$, and so
we may construct a larger group $\wmc$ the semi-direct product
$$\wmc:= \mc\sdp \Si$$
generated by the mapping class group and the sign changes.
Note that $\wmc$ actually acts on the set of conjugacy classes of 
triples of 
generators of the binary icosahedral group, and therefore also on 
the set $S$ of invariants of generating triples.  

\begin{defn}
Two seven-tuples $\bm,\bm'\in S$ are {\em geometrically equivalent} if
they are in the same orbit of the group 
$\wmc$.
%generated by the mapping class group $\mc$ and the group $\Si$ 
%of even sign changes.
\end{defn}

A key fact that we will use is:

\begin{lem} \label{lem: geq implies eq}
If two solutions of $\text{P}_{\text{VI}}$ have geometrically equivalent 
linear monodromy data in $S$, then the solutions are equivalent.
\end{lem}
\pf
First 
note that if two solutions have the same data $\bm$ then they are related
by a translation in $W_a(D^-_4)$.

Also recall that 
if two solutions have monodromy data related by  the free subgroup 
$\F_2\subset \mc$ then they are equivalent, since
the $\F_2$ action
corresponds to the branching of a single solution.

Thus it is sufficient, for 
each generator of $\wmc$, to find an Okamoto transformation
inducing the same action on the monodromy data, at least 
up to the action of $\F_2$.
To avoid confusion first note that there are two reasons that 
$W_a(F_4)$ does not in fact 
 act on the monodromy data:

First it is straightforward to check that $W_a(F_4)$
does not even act on the local monodromy 
data $(m_1,\ldots m_4)$; even the
subgroup $W_a(D^-_4)$ does not act here. This can be
easily rectified by working with the data $(\theta_i,m_{ij})$ instead.
Thus $W_a(D^-_4)$ acts on $\{(\theta_i,m_{ij})\}$  
(acting 
trivially on the quadratic
functions $m_{ij}$ by \cite{IIS}).

Secondly one still does not get an action of $W_a(F_4)$ on
$\{(\theta_i,m_{ij})\}$, since the $W_a(F_4)$ action on the systems 
\eqref{eq: lin syst} moves the pole positions, and so one obtains
representations of the fundamental group of different punctured spheres.
Although for each fixed $t$ only six four-punctured spheres arise (the $S_3$
orbit of $t$) one cannot
just add in $S_3$; one inevitably becomes involved with an infinite subgroup 
of $\mc$. (The essential reason for this is that the sequence 
\eqref{eq: mc to s4} does not split; this is implied by the non-splitting of 
\eqref{eq: modular} below, which follows from the Kurosh subgroup theorem.)
However as we will explain below, there is a well-defined action on the set of
$\F_2$ orbits in $\{(\theta_i,m_{ij})\}$ 
which is sufficient for us.

There are two steps: first the generators of $\mc$ will be related, up to the
$\F_2$ action, to certain automorphisms of $\IP^1$. Then the action
of these automorphisms on the systems \eqref{eq: lin syst} 
will be identified with Okamoto transformations. (The
action of the signs will be dealt with at the end.)

Suppose %$t\in B:=\IP^1\setminus\{0,1,\infty\}$ 
$f(z)$ is a M\"obius
transformation such that 
$f(\{0,t,1,\infty\})= \{0,t,1,\infty\}$
for some $t\in B:=\IP^1\setminus\{0,1,\infty\}$.
Thus $f$ represents an element of the mapping class group
$$\mc^t:= \pi_0(\Diff(S^2,\{0,t,1,\infty\}))$$
of the sphere with the points $\{0,t,1,\infty\}$ marked.
For example we will be interested in the cases:
$$
f_1=\frac{t-z}{t-1},\qquad
f_2=\frac{z}{t},\qquad
f_3=\frac{z}{z-1},$$
which represent elements of the groups $\mc^t$ when $t=2,-1,2$ respectively.
Now to identify these mapping class groups $\mc^t$ with a standard one, say
with $\mc=\mc^{1/2}$ we need to choose diffeomorphisms identifying these
four-pointed spheres; different choices of diffeomorphisms correspond to
conjugating by elements of $\F_2\subset \mc$.
Thus each $f_i$ leads to a well-defined element of $\mc/\F_2\cong S_4$.
These elements of $S_4$ are just the induced permutations of
$\{0,t,1,\infty\}$, i.e.
$$ f_i \text{ corresponds to the involution } (i,i+1)\in S_4.$$
Thus, up to $\F_2$,  the action of $f_i$ is given by the generator
$\om_i$ of $\mc$, for $i=1,2,3$.

Now we will identify the action, by pullback, of these M\"obius transformations
$f_i$ on the systems \eqref{eq: lin syst} in terms of Okamoto transformations.
This is straightforward since the Okamoto transformations are determined by
their action on $V=\{(\theta_i)\}$. 
We  find $f_1$ corresponds to the transformation
$x^3$ of \cite{OkaPVI} p.361 (with $1/(t-1)$ corrected to $t/(t-1)$), that
$f_2$ corresponds to $x^{313}:=x^3\circ x^1\circ x^3$
and $f_3$ corresponds to $s\circ x^{313}$ where $s\in W_a(D^-_4)$ is the
element acting on the $\theta$-parameters as the permutation $(14)(23)$. 
(One could check these directly or use the 
universality of Okamoto's action---that
all such M\"obius transformations will lead to Okamoto transformations.)

Finally we will obtain a generator of the sign changes (others being obtained
under the action of $S_4$).
Note that the transformation $x^2$ of \cite{OkaPVI} p.361 maps $\theta$
to $(\theta_4-1,\theta_2,\theta_3,\theta_1+1)$, so that
$T:=s\circ x^{313}\circ x^2$ maps $\theta$ to 
$(\theta_1+1,\theta_2,\theta_3,\theta_4-1)$.
We claim $T$ acts on the monodromy data just by negating 
$M_1$ and $M_4$ and fixing $M_2,M_3$. 
To see this we identify $T$  as a composition $T=\tau\circ R$.
Here $R$ 
is the rational gauge transformation (elementary
Schlesinger transformation \cite{JM81}) increasing by $1$ the first 
eigenvalue $\theta_1/2$ of the residue $A_1$ at zero, and
decreasing by $1$ the first eigenvalue $\theta_4/2$
of the residue $A_4$ at infinity. 
The resulting system has exactly the same monodromy data but is no longer 
in $\lsl_2$, so we apply the operation $\tau$  twisting
by the flat line bundle $-\frac{dz}{2z}$ to get an
$\lsl_2$  system with parameters 
$(\theta_1+1,\theta_2,\theta_3,\theta_4-1)$.  
The twisting gives the sign changes.
\epf

This leads to the main result:

\begin{thm} \label{thm: exactly 52}
There are exactly $52$ inequivalent icosahedral 
solutions to $\text{P}_{\text{VI}}$.
\end{thm}
\pf
First we compute, using Maple, the orbits of $\wmc$ in $S$.
We find there are exactly $52$ orbits.
Thus by Lemma \ref{lem: geq implies eq}
 there are at most $52$ inequivalent icosahedral solutions
to $\text{P}_{\text{VI}}$.
Combining this with Proposition \ref{prop: at least 52} yields result.
\epf

Some properties of these solutions will be listed in the next section.

\begin{rmk}\label{rmk: lin alg solns}
One may be interested in classifying the systems \eqref{eq: lin syst}
having bases of algebraic solutions.
The natural equivalence relation to use is geometric equivalence (of the
correpsonding monodromy representations) since this preserves the (projective)
monodromy group.
Theorem \ref{thm: exactly 52} implies one will get the same classification as
that appearing here (for the icosahedral representations).
There is still the thorny question of determining precisely which pole 
positions 
are possible for systems having given exponents and monodromy representation.
But this is determined in a straightforward way if we know the corresponding
\PVI\ solution (namely the \PVI\ solution explicitly determines an
isomonodromic family of systems and one examines when there are poles in the
matrix entries of these systems).
See also remark \ref{rmk: lame}.
\end{rmk}

\begin{rmk}\ 

\noindent
1) In general parameter equivalence is strictly weaker than equivalence,
even when restricted to algebraic solutions 
(cf. e.g. \cite{Hit-Poncelet}). Adding other invariants (such as the
number of branches, the genus, the nonlinear monodromy group size) 
does often distinguish algebraic solutions, but not always:
For example \eqref{eq: dih soln} and \eqref{eq: hitoct} below are both
solutions for $(\theta_1,\theta_2,\theta_3,\theta_4)=(1,1,1,3)/6$ and have
identical invariants, but are inequivalent (since they are
not related by the
subgroup of $W_a(F_4)$ stabilising $(1,1,1,3)/6$; each of the three
reflections $\si_{03}, \si_{04}$ and $s_1s_2s_1$ of \cite{Masuda-PVI, NY-so8}
fix both solutions).

\noindent
2) In general geometric equivalence is strictly stronger than equivalence;
one needs to add the action of the rest of $W_a(D_4^-)$. This is slightly
subtle since, as mentioned above, $W_a(D_4^-)$ {does not act} on the set
$\Hom(\pi_1(S^*),G)/G$ 
of conjugacy classes of
fundamental group representations---one
needs to either use a covering such as 
$\{(\theta_i,m_{ij})\}$ (or a suitable finite intermediate cover) or a quotient
(such as the coefficients of the Fricke relation and the $m_{ij}$, on which 
$W_a(D_4^-)$ acts trivially).
\end{rmk}

\begin{rmk}
Usually \cite{DubMaz00, Iwas-modular, k2p} the $\F_2$ action on the 
monodromy data
is described in terms of the three-string (Artin/planar) braid groups
$$1\to P_3\to B_3 \to S_3\to 1$$
whose centres act trivially. Quotienting this sequence by the centres
$Z(B_3)=Z(P_3)\cong\IZ$ one obtains
\beq\label{eq: modular}
1\to\F_2 \to \PSL_2(\IZ)\to S_3 \to 1
\eeq
so $\F_2$ is identified with the level-two subgroup $\Ga(2)$.
However to obtain the full symmetries of $\text{P}_{\text{VI}}$ one needs to extend $\F_2$ by
$S_4$, as in \eqref{eq: mc to s4}, and not just by $S_3$; this is why we used
$\mc$. 
This just corresponds to pulling back \eqref{eq: modular} along the natural
map $S_4\to S_3$ with kernel $K_4$ the Klein four group.
Indeed there is a map $\mc\to\PSL_2(\IZ)$ with kernel $K_4$
(arising from the 
cross-ratio $(\IP^1)^4\setminus\text{diags}\to B$), 
and in fact $\mc\cong\PSL_2(\IZ)\sdp K_4$ (cf. \cite{Birman} p.206).
\end{rmk}

\end{section}

\begin{section}{Properties of the $52$ solutions}\label{sn: tables}

In this section we will list some of the properties of the $52$ icosahedral
solutions to $\text{P}_{\text{VI}}$.
See tables $1$ and $2$.
The columns of table $1$ are defined as follows:

{\Small
\begin{table}
\begin{center}

\begin{tabular}[h]{|c|c|c|c|c|c|c|c|c|c| }
\hline
  & \text{Degree}  &  \text{Genus} & \text{Walls} & \text{$A_5$ Type} & Alcove Point
  & n &
  Good? &  Group (Size) & Partitions \\ \hline
1 & 1 & 0 & 1 & $ a\,b\,c$ & 31, 19, 11, 1 & 192 & $\circ$ & $1$ & $ $\\ \hline
2 & 1 & 0 & 1 & $ a\,b\,d$ & 37, 17, 13, 7 & 192 & $\circ$ & $1$ & $ $\\ \hline
3 & 1 & 0 & 1 & $ a\,c\,d$ & 33, 21, 9, 3 & 192 & $\circ$ & $1$ & $ $\\ \hline
4 & 1 & 0 & 1 & $ b\,c\,d$ & 28, 16, 8, 4 & 192 & $\circ$ & $1$ & $ $\\ \hline
5 & 1 & 0 & 2 & $ b^2\,c$ & 26, 14, 6, 6 & 96 & $\circ$ & $1$ & $ $\\ \hline
6 & 1 & 0 & 2 & $ b^2\,d$ & 38, 18, 18, 2 & 96 & $\circ$ & $1$ & $ $\\ \hline
7 & 1 & 0 & 2 & $ b\,c^2$ & 22, 10, 10, 2 & 96 & $\circ$ & $1$ & $ $\\ \hline
8 & 1 & 0 & 2 & $ b\,d^2$ & 34, 14, 10, 10 & 96 & $\circ$ & $1$ & $ $\\ \hline
9 & 1 & 0 & 3 & $ c^3$ & 18, 6, 6, 6 & 32 & $\circ$ & $1$ & $ $\\ \hline
10 & 1 & 0 & 3 & $ d^3$ & 42, 18, 18, 6 & 32 & $\circ$ & $1$ & $ $\\ \hline
11 & 2 & 0 & 2 & $ b^2\,c^2$ & 42, 18, 10, 10 & 96 & $\times$ & $2$ & $1$, ${2}$\\ \hline
12 & 2 & 0 & 2 & $ b^2\,d^2$ & 50, 10, 6, 6 & 96 & $\times$ & $2$ & $1$, ${2}$\\ \hline
13 & 2 & 0 & 2 & $ c^2\,d^2$ & 42, 18, 6, 6 & 96 & $\times$ & $2$ & $1$, ${2}$\\ \hline
14 & 3 & 0 & 1 & $ b\,c^2\,d$ & 40, 16, 8, 8 & 288 & $\times$ & $S_{3}$ & ${3}$, ${2}$\\ \hline
15 & 3 & 0 & 1 & $ b\,c\,d^2$ & 40, 8, 4, 4 & 288 & $\times$ & $S_{3}$ & ${3}$, ${2}$\\ \hline
16 & 4 & 0 & 2 & $ a\,c^3$ & 33, 9, 9, 9 & 128 & $\circ$ & $A_{4}$ & ${3}$\\ \hline
17 & 4 & 0 & 2 & $ a\,d^3$ & 51, 3, 3, 3 & 128 & $\circ$ & $A_{4}$ & ${3}$\\ \hline
18 & 4 & 0 & 2 & $ c^3\,d$ & 30, 6, 6, 6 & 128 & $\circ$ & $A_{4}$ & ${3}$\\ \hline
19 & 4 & 0 & 2 & $ c\,d^3$ & 42, 6, 6, 6 & 128 & $\circ$ & $A_{4}$ & ${3}$\\ \hline
20 & 5 & 0 & 1 & $ b^2\,c\,d$ & 44, 12, 12, 4 & 480 & $\times$ & $S_{5}$ & ${2}^{2}$, ${2}\,{3}$\\ \hline
21 & 5 & 0 & 2 & $ c^2\,d^2$ & 36, 12, 0, 0 & 240 & $\times$ & $S_{5}$ & ${3}$, ${2}\,{3}$\\ \hline
22 & 6 & 0 & 1 & $ b\,c^2\,d$ & 34, 10, 2, 2 & 576 & $\circ$ & $S_{6}$ & ${5}$, ${2}\,{3}$\\ \hline
23 & 6 & 0 & 1 & $ b\,c\,d^2$ & 46, 14, 10, 2 & 576 & $\circ$ & $S_{6}$ & ${5}$, ${2}\,{3}$\\ \hline
24 & 8 & 0 & 1 & $ a\,c^2\,d$ & 39, 15, 3, 3 & 768 & $\times$ & $A_{8}$ & ${3}\,{5}$, ${2}^{2}\,{3}$\\ \hline
25 & 8 & 0 & 1 & $ a\,c\,d^2$ & 45, 9, 9, 3 & 768 & $\times$ & $A_{8}$ & ${3}\,{5}$, ${2}^{2}\,{3}$\\ \hline
26 & 9 & 1 & 2 & $ b\,c^3$ & 28, 4, 4, 4 & 288 & $\circ$ & $A_{9}$ & ${3}\,{5}$\\ \hline
27 & 9 & 1 & 2 & $ b\,d^3$ & 52, 8, 8, 4 & 288 & $\circ$ & $A_{9}$ & ${3}\,{5}$\\ \hline
28 & 10 & 0 & 2 & $ a^2\,c\,d$ & 48, 12, 6, 6 & 480 & $\times$ & ${2}^{7}\,{3}\,{5}$ & ${2}^{2}\,{3}^{2}$\\ \hline
29 & 10 & 0 & 2 & $ b^3\,c$ & 46, 14, 14, 6 & 320 & $\circ$ & $A_{10}$ & ${2}^{2}\,{5}$\\ \hline
30 & 10 & 0 & 2 & $ b^3\,d$ & 42, 2, 2, 2 & 320 & $\circ$ & $A_{10}$ & ${2}^{2}\,{5}$\\ \hline
31 & 10 & 0 & 3 & $ c^4$ & 24, 0, 0, 0 & 80 & $\circ$ & $A_{10}$ & ${3}\,{5}$\\ \hline
32 & 10 & 0 & 3 & $ d^4$ & 48, 0, 0, 0 & 80 & $\circ$ & $A_{10}$ & ${3}\,{5}$\\ \hline
33 & 12 & 0 & 0 & $ a\,b\,c\,d$ & 43, 11, 7, 1 & 2304 & $\times$ & $A_{12}$ & ${2}^{2}\,{3}^{2}$, ${2}^{2}\,{3}\,{5}$\\ \hline
34 & 12 & 1 & 1 & $ a\,b\,c^2$ & 37, 13, 5, 5 & 1152 & $\times$ & $A_{12}$ & ${3}^{2}\,{5}$, ${2}^{2}\,{3}\,{5}$\\ \hline
35 & 12 & 1 & 1 & $ a\,b\,d^2$ & 49, 5, 5, 1 & 1152 & $\times$ & $A_{12}$ & ${3}^{2}\,{5}$, ${2}^{2}\,{3}\,{5}$\\ \hline
36 & 12 & 1 & 1 & $ b^2\,c\,d$ & 38, 6, 6, 2 & 1152 & $\times$ & ${2}^{9}\,{3}^{2}\,{5}$ & ${2}^{2}\,{3}^{2}$, ${2}\,{5}^{2}$\\ \hline
37 & 15 & 1 & 2 & $ b^3\,c$ & 36, 4, 4, 4 & 480 & $\times$ & $A_{15}$ & ${2}^{2}\,{3}^{2}\,{5}$\\ \hline
38 & 15 & 1 & 2 & $ b^3\,d$ & 48, 8, 8, 8 & 480 & $\times$ & $A_{15}$ & ${2}^{2}\,{3}^{2}\,{5}$\\ \hline
39 & 15 & 1 & 2 & $ b^2\,c^2$ & 32, 8, 0, 0 & 720 & $\times$ & $S_{15}$ & ${3}^{2}\,{5}$, ${2}\,{3}\,{5}^{2}$\\ \hline
40 & 15 & 1 & 2 & $ b^2\,d^2$ & 44, 4, 0, 0 & 720 & $\times$ & $S_{15}$ & ${3}^{2}\,{5}$, ${2}\,{3}\,{5}^{2}$\\ \hline
41 & 18 & 1 & 3 & $ b^4$ & 40, 0, 0, 0 & 144 & $\circ$ & ${2}^{14}\,{3}^{4}\,{5}\,{7}$ & ${3}^{2}\,{5}^{2}$\\ \hline
42 & 20 & 1 & 1 & $ a\,b^2\,c$ & 41, 9, 9, 1 & 1920 & $\times$ & $A_{20}$ & ${2}^{4}\,{3}^{2}\,{5}$, ${2}^{2}\,{3}^{2}\,{5}^{2}$\\ \hline
43 & 20 & 1 & 1 & $ a\,b^2\,d$ & 47, 7, 3, 3 & 1920 & $\times$ & $A_{20}$ & ${2}^{4}\,{3}^{2}\,{5}$, ${2}^{2}\,{3}^{2}\,{5}^{2}$\\ \hline
44 & 20 & 1 & 3 & $ a^2\,c^2$ & 42, 18, 0, 0 & 480 & $\times$ & ${2}^{17}\,{3}^{4}\,{5}^{2}\,{7}$ & ${3}^{2}\,{5}^{2}$, ${2}^{2}\,{3}^{2}\,{5}^{2}$\\ \hline
45 & 20 & 1 & 3 & $ a^2\,d^2$ & 54, 6, 0, 0 & 480 & $\times$ & ${2}^{17}\,{3}^{4}\,{5}^{2}\,{7}$ & ${3}^{2}\,{5}^{2}$, ${2}^{2}\,{3}^{2}\,{5}^{2}$\\ \hline
46 & 24 & 1 & 2 & $ a\,b^3$ & 45, 5, 5, 5 & 768 & $\times$ & ${2}^{20}\,{3}^{5}\,{5}^{2}\,{7}\,{11}$ & ${2}^{4}\,{3}^{2}\,{5}^{2}$\\ \hline
47 & 30 & 2 & 2 & $ a^2\,b\,c$ & 46, 14, 4, 4 & 1440 & $\times$ & ${2}^{24}\,{3}^{6}\,{5}^{3}\,{7}^{2}\,{11}\,{13}$ & ${2}^{2}\,{3}^{2}\,{5}^{4}$, ${2}^{4}\,{3}^{4}\,{5}^{2}$\\ \hline
48 & 30 & 2 & 2 & $ a^2\,b\,d$ & 52, 8, 2, 2 & 1440 & $\times$ & ${2}^{24}\,{3}^{6}\,{5}^{3}\,{7}^{2}\,{11}\,{13}$ & ${2}^{2}\,{3}^{2}\,{5}^{4}$, ${2}^{4}\,{3}^{4}\,{5}^{2}$\\ \hline
49 & 36 & 3 & 3 & $ a^2\,b^2$ & 50, 10, 0, 0 & 864 & $\times$ & ${2}^{23}\,{3}^{4}\,{5}\,{7}$ & ${3}^{4}\,{5}^{4}$, ${2}^{2}\,{3}^{4}\,{5}^{4}$\\ \hline
50 & 40 & 3 & 3 & $ a^3\,c$ & 51, 9, 9, 9 & 320 & $\times$ & ${2}^{25}\,{3}^{4}\,{5}^{2}\,{7}$ & ${2}^{4}\,{3}^{4}\,{5}^{4}$\\ \hline
51 & 40 & 3 & 3 & $ a^3\,d$ & 57, 3, 3, 3 & 320 & $\times$ & ${2}^{25}\,{3}^{4}\,{5}^{2}\,{7}$ & ${2}^{4}\,{3}^{4}\,{5}^{4}$\\ \hline
52 & 72 & 7 & 3 & $ a^3\,b$ & 55, 5, 5, 5 & 576 & $\times$ & ${2}^{32}\,{3}^{4}\,{5}\,{7}$ & ${2}^{4}\,{3}^{8}\,{5}^{8}$\\ \hline

\hline
\end{tabular}

\vspace{0.2cm}
\caption{Properties of the 52 icosahedral solutions.}
\label{solution table}
\end{center}\end{table}
%\end{landscape}
}

\noindent$\bullet$\ 
The degree is the number of branches that the solution has. (Recall
solutions branch at $t\in\{0,1,\infty\}$.)

\noindent$\bullet$\ 
The genus is the genus of the algebraic curve on which the function $y(t)$
becomes single-valued. This is computed using the Riemann--Hurwitz formula
from the permutation representation of the cover.

\noindent$\bullet$\ 
The column labelled `Walls' lists the number of affine $F_4$ hyperplanes that
the parameters of the solutions lie on. Since the Okamoto transformations
reflect in these hyperplanes, this number is an invariant. 

\noindent$\bullet$\ 
The $A_5$ type of the solution is defined as follows. Recall the icosahedral
rotation group $A_5$ has precisely five conjugacy classes. We label the four
non-trivial classes, the rotations by 
$1/2, 1/3, 1/5, 2/5$ of a turn, by the letters $a,b,c,d$ respectively.
Thus given a four-tuple $M_1,M_2,M_3,M_4$ of elements of the binary
icosahedral group $\Ga$ 
we are listing the set of conjugacy classes of their image
in $A_5=\Ga/\pm$. (If there are only three classes listed, that means that the
fourth class is the trivial class.)
This set is an invariant of the icosahedral solution
although it does not determine the equivalence class of the solution (compare
e.g. rows 12 and 40).

\noindent$\bullet$\ 
The alcove point is the value of (sixty times) the unique four-tuple of
equivalent parameters $(\theta_1,\theta_2,\theta_3,\theta_4)$ which lies in the
closure of the standard alcove \eqref{eq: std alc}. We scale by $60$ simply to
clear the denominators. This is the `parameter equivalence class'.

\noindent$\bullet$\ 
The value $n$ in the next column is the number of $7$-tuples $\bm\in S$ 
corresponding to that equivalence class. 
Thus each $n$ is divisible by the corresponding degree
and the sum of all the $n$'s is $26688$. 

\noindent$\bullet$\ 
Let $k$ be the degree of one of the solutions.
We will say the solution is `good'
if it has a representative (amongst the $n/k$ coming from $S$)
for which Jimbo's formula (cf. \cite{Jimbo82} and \cite{k2p} Theorem 4) may be
applied to give the asymptotics at $t=0$ on every branch.
A cross ($\times$) means the solution has such a good representative and that
we can in principle apply the procedure of \cite{k2p} to find the solution
explicitly. Even if a solution is not good ($\circ$) there may well 
be other ways to identify the solution (see below).

\noindent$\bullet$\ 
The column `Group (Size)' lists the nonlinear monodromy group of the solution
or at least the size of this group;
this is the group generated by the permutations of the branches of the
solution curve as $t$ goes around $0,1,\infty$ 
and so naturally appears as a subgroup
of the symmetric group on $k$ letters. (In other words it 
is the monodromy group of the solution curve, expressed as a branched
cover of $\IP^1\setminus\{0,1,\infty\}$, or equivalently the Galois group of
this cover over the base field $\IC(t)$.)

\noindent$\bullet$\ 
The final column `Partitions' lists the set of 
conjugacy classes of the three generators of the nonlinear monodromy group.
These are conjugacy classes in $\Sym_k$, where $k$ is the degree of the
solution, and so are written as partitions of $k$ (representing the cycle
lengths of the permutations). If there are less than $3$
partitions listed (separated by commas) the {\em last} partition is repeated.
E.g. in row $20$ the degree is $5$ and the three partitions indicated are
$1+2+2, 2+3, 2+3$ (repeating the last one). The set of these partitions is an
invariant.

\begin{rmk}
Observe that there are a number of consecutive rows of table $1$
which look the same except
for having a different alcove point and having $c,d$ swapped in their $A_5$
type. We will refer to these as {\em sibling}
solutions. They occur because $A_5$
has a non-trivial outer automorphism (from the extension
$1\to A_5\to S_5\to \IZ/2\to 1$, i.e. by conjugating an element of $A_5$ by an
odd permutation) which swaps the conjugacy classes $c,d$ and preserves the
others. The action of this outer automorphism on the additive side, 
on the solutions themselves,
remains mysterious.
\end{rmk}

\begin{rmk}
Observe that the first ten entries of table $1$ correspond to solutions with
only one branch. Looking at the $A_5$ type we see that in each case one of the
local (linear) monodromies is (projectively) trivial. Thus these ten correspond
to the
list of equivalence  classes of {\em pairs} of generators of the
icosahedral group, i.e. to hypergeometric equations with icosahedral
monodromy, i.e. to the ten icosahedral entries on Schwarz's list
\cite{Schwarz}. 
(Replacing
$a$ by $1/2$, $b$ by $1/3$ etc. gives the bijection with this part of 
Schwarz's list.)
As solutions of $\text{P}_{\text{VI}}$ these are all equivalent to a constant solution
(and also equivalent to the solution $y=t$, with parameters as listed 
in table $2$). 
Thus it is tempting to view 
Okamoto's $W_a(F_4)$-action
as
the natural extension (to linear systems of the form \eqref{eq: lin syst}) 
of the equivalence relation used by Schwarz.
But this is not quite right as Okamoto's action does not preserve the linear
monodromy group (cf. remark \ref{rmk: lin alg solns}).
Rather the $W_a(F_4)$ action is the natural analogue for the nonlinear
\PVI\ equation
of Schwarz's equivalence relation 
(and this was how Okamoto was thinking of it in 
\cite{OkaPVI}); it does indeed preserve algebraicity of \PVI\ solutions. 
The extension of the 
further step taken by Klein, of obtaining the entries of
Schwarz's list as pullbacks
along a rational map
from the `basic Schwarz list', is also possible (indeed Klein proves this, at
least for the corresponding Fuchsian equations).
However to construct such rational pullbacks explicitly, before the system
upstairs has
been found explicitly, is a difficult problem (moreover 
the $\text{P}_{\text{VI}}$ solution is equivalent to constructing a complete
algebraic family of such covers as the four pole positions move). 
Such a procedure has been described (modulo the difficulty mentioned) 
independently 
by Doran \cite{Chuck1} and Kitaev \cite{Kitaev-dessins}.
(Kitaev's paper \cite{Kitaev-dessins} also contains
some new explicit 
examples---i.e. not equivalent to solutions previously
constructed by other means---see 2) of 
Remark \ref{rmk: known solutions} below.)
\end{rmk}

\begin{rmk}
The next solutions ($\leqslant 4$ branches) 
are simple deformations of known solutions, i.e. the same solution
just with different parameters, as follows.

The solutions with two branches ($11,12$ and $13$) 
are equivalent to the solution
$y=\pm\sqrt{t}$. (One first observes that if 
$\theta_2=\theta_3$ and $\theta_1+\theta_4=1$ then 
this is indeed a solution. Then one checks each of solutions $11,12,13$ has a
representative with such parameters, as listed in table $2$. 
Finally one uses Jimbo's formula to see
that a leading term of these solutions matches that of $y=\pm\sqrt{t}$.)

Similarly the three-branch tetrahedral solution
\begin{eqn}\label{eq: hittet}
y=\frac{(s-1)(s+2)}{s(s+1)},\qquad
t=\frac{(s-1)^2(s+2)}{(s+1)^2(s-2)}
\end{eqn}
\!\!on p.592 of \cite{Hit-Octa}
is actually a 
solution on the whole line 
$\theta_1/2=\theta_2=\theta_3,\theta_4=\frac{2}{3}$, 
amongst other possibilities.
(Note $\be$ should be  $-2/9$, not $-1/18$ in op. cit.)
Solutions $14$ and $15$ have representatives on this line (see table $2$).
Their leading terms given by Jimbo's formula (on the two branches where it may
be applied)
are $\pm i\sqrt{3}t^{1/2}$,
which match the Puiseux expansion of \eqref{eq: hittet}, so \eqref{eq: hittet}
gives both icosahedral solutions $14$ and $15$. (Using the $\text{P}_{\text{VI}}$ equation
the leading terms determine the whole Puiseux expansion and thus the entire 
solution.)

Next the four-branch dihedral solution in section 6.1 of \cite{Hit-Poncelet}:
\begin{equation} \label{eq: dih soln}
y=\frac{s^2(s+2)}{s^2+s+1},\qquad
t=\frac{s^3(s+2)}{2s+1}
\end{equation}
is a solution if $\theta_1=\theta_2=\theta_3, \theta_4=1/2$.
As above this gives the icosahedral solutions $16$ and $17$, 
with the parameters indicated in table $2$.

Finally the four-branch octahedral solution
\begin{eqn}\label{eq: hitoct}
y=\frac{(s-1)^2}{s(s-2)},\qquad
t=\frac{(s+1)(s-1)^3}{s^3(s-2)}
\end{eqn}
on p.588 of \cite{Hit-Octa}
is a solution if $\theta_1=\theta_2=\theta_3$ and either
$$\theta_4=1-3\theta_1\quad\text{or}\quad\theta_4=1+3\theta_1.$$ 
(The implicit version of this in \cite{Hit-Octa} should read: 
$3\,{y}^{4}- \left( 4\,t+4 \right) {y}^{3}+6\,t{y}^{2}-{t}^{2}=0.$)

This gives icosahedral solution $18$, with parameters as in table $2$.
Solution $19$ is slightly more elusive and looks not to have a representative 
in either family. 
However it is {\em equivalent} to a member of the second family:
Take the solution \eqref{eq: hitoct}, with 
$(\theta_1,\theta_2,\theta_3,\theta_4)=(2,2,2,11)/5$.
Then apply the sequence of Okamoto transformations
$$s_1(s_2s_0s_3s_4)^2s_2$$
in the notation of \cite{NY-so8} (we act on the left, so we do the right-most 
$s_2$ first). This is the transformation which reduces
$\theta_4$ by $2$ so yields an explicit solution with
parameters $(2,2,2,1)/5$, which may be parameterised as follows: 
$$y={\frac {7+22\,s+7\,{s}^{2}}{ 8\,\left( 1+s+{s}^{2} \right) s \left( 
s+2 \right) }},\qquad
t=
{\frac {1+2\,s}{{s}^{3} \left( s+2 \right) }}.$$
The corresponding Puiseux expansion at zero has a leading term 
$7\,\, 2^{1/3}t^{2/3}/16$ which matches the leading term given by Jimbo's
formula for entry $19$ of table $2$, and so this is icosahedral solution $19$.

(We remark that there is thus a problem with the nomenclature for the
expressions \eqref{eq: hittet}-\eqref{eq: hitoct}, they are as much
{\em icosahedral} as they are {\em tetrahedral} etc. 
This is on top of the fact that 
Hitchin's octahedral solution \eqref{eq: hitoct} 
is equivalent to a solution found independently by Dubrovin
\cite{Dub95long} starting from the tetrahedral reflection group, and similarly
for \eqref{eq: hittet} and the octahedral reflection group, cf. \cite{k2p} 
Remark 14.)
\end{rmk}

\begin{rmk}
The three Dubrovin--Mazzocco icosahedral solutions  
\cite{Dub95long, DubMaz00} 
are equivalent
to the solutions on rows $31,32$ ($10$ branch siblings) and $41$ ($18$ branches,
genus one).
To prove they are equivalent one can show that the unipotent monodromy data
used in  \cite{DubMaz00} 
may be mapped by an Okamoto transformation to a triple of
generators of $\Ga$ (cf. \cite{k2p} Remark 14).
(It is sufficient to use the affine $D_4$ group, which acts trivially on the
quadratic functions $m_{ij}$ of the monodromy data \cite{IIS}.)
Alternatively a simpler but less direct way to see this is to observe
that the icosahedral solutions $31,32,41$ here 
have equivalent parameters to those
of \cite{DubMaz00}. Then appeal to the classification of \cite{DubMaz00}  
of all such finite branching solutions. 
Observe also that apart from this $18$ branch solution, 
all icosahedral solutions with more than 
$10$ branches are good, and so their Puiseux expansions at $0$ may be computed
using Jimbo's formula. Solution $41$ appears in Theorem C
of the introduction; Jimbo's formula yields the leading term on $16$ of the
$18$ branches at zero and for the other two (where the solution does not in
fact branch) we used Okamoto transformations to transfer 
the leading terms given by Dubrovin and
Mazzocco \cite{DubMaz00} (in fact one needs the first $2$ terms in the Taylor
expansion in order for $\text{P}_{\text{VI}}$ to determine the series uniquely, and to compute
these terms one needs the first $3$ terms of the corresponding branches of
Dubrovin--Mazzocco's solution, but these are easily found from the given 
leading terms in \cite{DubMaz00}).
The equivalent parameters used in \cite{DubMaz00} were
$\theta_1=\theta_2=\theta_3=0, \theta_4=-2/3$ and by using Okamoto
transformations it is straightforward to convert the parameterisation of
Theorem C into a solution for these parameters. The result is:

$$y=\frac{1}{2}+{\tiny\frac{
\left(\begin{matrix}
128\,{s}^{18}-2496\,{s}^{17}+19728\,{s}^{16}+4605216\,{s}
^{15}-53030400\,{s}^{14}+229874976\,{s}^{13}-600089472\,{s}^{12}+\\
968994816\,{s}^{11}-
823777848\,{s}^{10}-88169600\,{s}^{9}+1204313064\,
{s}^{8}-1658437668\,{s}^{7}+1282505784\,{s}^{6}-\\
632776452\,{s}^{5}+
199216125\,{s}^{4}-36900918\,{s}^{3}+3168636\,{s}^{2}+134172\,s-38416
\end{matrix}\right)}
{6u\left(\begin{matrix}
5776\,{s}^{15}-85440\,{s}^{14}+482880\,{s}^{13}-1490080\,{s}^
{12}+13986240\,{s}^{11}-58604928\,{s}^{10}+\\133381480\,{s}^{9}-
186525360\,{s}^{8}+
162484560\,{s}^{7}-80442380\,{s}^{6}+11088528\,{s}^
{5}+\\12426960\,{s}^{4}-9203395\,{s}^{3}+3037020\,{s}^{2}-496860\,s+
33124 
\end{matrix}\right)}}
$$
with $t,u,s$ exactly as in Theorem C.
\end{rmk}

\begin{rmk}\label{rmk: known solutions}
Even if a solution is not `good' it may well be accessible: 

1) As already discussed, the smaller
solutions (1-4 branches) 
are simple deformations of known solutions.

2) Page $12$ of A. Kitaev's paper \cite{Kitaev-dessins}
contains an explicit formula for the solution
on row $26$ of table $1$, the smallest genus one solution. Presumably the
sibling solution (row $27$) can be obtained similarly;
in any case we will obtain it with our methods in section \ref{sn: moreex}. 
(Also \cite{Kitaev-dessins} (3.3) p.24 corresponds to row 21.) 

3) The Dubrovin--Mazzocco icosahedral solutions are not good in the sense of
   table $1$, but were found by
   adapting Jimbo's formula to their situation. 
   A different adaptation will be  made at the end of section \ref{sn: moreex}
   to find the asymptotics of the outstanding solutions.
\end{rmk}

\begin{rmk}\label{rmk: lame}
The three largest solutions, rows $50,51$ (genus three, $40$ 
branch siblings) and row
$52$ (genus seven, $72$ branches), are related to Lam\'e equations (certain 
second order Fuchsian ordinary differential equations 
having four singularities on $\IP^1$, and no
apparent singularities).
Namely there are Lam\'e equations having these (projective) monodromy
representations, given explicitly in the paper \cite{B-vdW} of 
Beukers and van der Waall.
Converting their equations into Fuchsian systems will give initial
conditions for these three $\text{P}_{\text{VI}}$ solutions, 
although evolving $\text{P}_{\text{VI}}$ to get a closed
form for these solutions is somewhat daunting.
(The leading terms in the Puiseux expansions given by Jimbo's formula
seem to give more information however.)
Note that the corresponding isomonodromic deformation will not be within the
space of Lam\'e equations---the deformed equations will have an apparent 
singularity.
Said differently, we may follow R. Fuchs and
think of \PVI\ as controlling isomonodromic deformations of rank two
Fuchsian equations (rather than systems), having four non-apparent
singularities at $z=0,t,1,\infty$ and an apparent singularity at $z=y$
(cf. e.g. \cite{NY-so8} for the formulae).
If we choose $t$ such that $y=0,1,t$ or $\infty$ then this equation will have
only the four non-apparent
singularities, so will be a Heun equation (with finite monodromy) and for 
the last three solutions this Heun equation will be of Lam\'e type.
All Heun/Lam\'e equations with finite monodromy should arise in this way.
\end{rmk}

%\begin{rmk}
As an example let us list the leading term at zero of the asymptotic
expansion of each branch of solution $52$. (We consider the representative with
$\theta_1=\theta_2=\theta_3=1/2, \theta_4=2/3$.)
The leading terms are each of the form $c\times t^{1-\si}$ 
where the coefficients $c$ are 
given by Jimbo's formula (\cite{Jimbo82}, \cite{k2p} Theorem 4).
To express the coefficients as algebraic numbers we raise them (or their real/imaginary
parts) to sufficiently high powers until they become rational and then look at
the continued fraction expansions. 
From table $1$, the $72$ branches over zero 
are partitioned into 4 two-cycles, 8 three-cycles and 8 five-cycles.
The values of $\sigma$ and $c$ for one branch of
each cycle are as follows:

{
%\begin{align}
%\si&=1/2:\ 
%&&\frac{1}{2}\,\sqrt {3\,\sqrt {5}+6\,i}\quad
%&&\frac{1}{2}\,\sqrt {3\,\sqrt {5}-6\,i}\quad
%&&\frac{1}{2}\,\sqrt {-3\,\sqrt {5}-6\,i}\quad
%&&\frac{1}{2}\,\sqrt {-3\,\sqrt {5}+6\,i}\notag
%\si&=1/3:\ 
%&&{2}^{3/5} {e^{7\pi\,i/10 }}\,\al\quad 
%&&{2}^{3/5} {e^{-3\pi\,i/10 }}\,\al\quad
%&&{2}^{3/5} {e^{3\pi\,i/10 }}\,\be\quad 
%&&{2}^{3/5} {e^{-7\pi\,i/10 }}\,\be \notag \\
%\si&=2/3:\ 
%&&\frac{2}{3}\,{2}^{{\frac {4}{15}}}{e^{7\pi\,i/10 }}\,\al\quad
%&&\frac{2}{3}\,{2}^{{\frac {4}{15}}}{e^{-3\pi\,i/10 }}\,\al\quad
%&&\frac{2}{3}\,{2}^{{\frac {4}{15}}}{e^{3\pi\,i/10 }}\,\be\quad
%&&\frac{2}{3}\,{2}^{{\frac {4}{15}}}{e^{-7\pi\,i/10 }}\,\be\notag 
%\end{align}
$$
\si=\frac{1}{2}:\ 
\frac{1}{2}\left(\pm 3\,\sqrt {5}\pm 6\,i\right)^{\frac{1}{2}},
\quad
\si=\frac{1}{3}:\ 
\pm \left(6\sqrt{3}\pm 2i\sqrt{5}\right)^{\frac{1}{3}},
\quad
\si=\frac{2}{3}:\ 
\pm \frac{2}{3}\left(3\sqrt{3}\pm i\sqrt{5}\right)^{\frac{1}{3}},
$$
$$
\si=\frac{1}{5}:\ 
\pm 2\,i\,{6}^{3/5},\quad
\si=\frac{2}{5}:\ 
\pm{\frac {9}{7}}\,i\,{6}^{1/5},\quad
\si=\frac{3}{5}:\ 
\pm{\frac {4}{13}}\,i\,{6}^{4/5},\quad
\si=\frac{4}{5}:\ 
\pm{\frac {12}{19}}\,i\,{6}^{2/5}.
$$
}
%\noindent{\!\!where} 
%$\al=\left( 79178+26970\,i\sqrt {15} \right) ^{1/30}$ and 
%$\beta=\left(79178-26970\,i\sqrt {15} \right) ^{1/30}$.
To obtain the other leading terms for each cycle, one just multiplies by all 
the
$k$-th roots of unity, where $k$ is the cycle length.
(Here $\si$ is the same for all branches of each cycle and
the cycle length is equal to the denominator of $\si$.)
It is still a challenge to write down the polynomial $F(t,y)=0$,
of degree $72$ in $y$
defining the solution curve and having these leading terms in its 
Puiseux expansions over $t=0$.
The curve itself is determined abstractly 
by its permutation representation as a cover of
$\IP^1\setminus\{0,1,\infty\}$. This can easily be computed applying the
operations $\omega^2_i$ to the seven-tuples of monodromy data.
The corresponding permutations of the $72$ branches thus obtained, going
around the loops $w_1,w_2$ of section \ref{sn: geom}, are, respectively:

{\Small
\begin{align*}
(1\, 2)(3\, 4)(5\, 6)(7\, 8)&
( 9\, 10\, 11)(12\, 13\, 14)(15\, 16\, 17)(18\, 19\, 20)
(21\, 22\, 23)(24\, 25\, 26)(27\, 28\, 29)(30\, 31\, 32)\\ \notag
&(33\ldots 37)(38\ldots 42)
(43\ldots 47)(48\ldots 52)
(53\ldots 57)(58\ldots 62)
(63\ldots 67)(68\ldots 72)\notag
\end{align*}
\begin{align*}
(1\, 18\, &39\, 49\, 21)(2\, 15\, 33\, 45\, 24)(3\, 22\, 59\, 68\, 20)(4\, 25\, 56\, 65\, 17)
(5\, 27\, 50\, 38\, 11)(6\, 30\, 46\, 37\, 14)(7\, 9\, 69\, 58\, 29)(8\, 12\,
  66\, 55\, 32)\\ \notag
&(10\, 34\, 64)(13\, 40\, 72)(16\, 70\, 42)(19\, 67\, 36)(23\, 47\, 54)(26\, 51\, 62)
(28\, 57\, 44)(31\, 60\, 48)(35\, 52)(41\, 43)(53\, 71)(61\, 63)\notag
\end{align*}
}
One easily computes, using Riemann--Hurwitz, that this represents a genus $7$
Belyi curve, and wonders whether this curve is remarkable for any other reasons.
%\end{rmk}

To end this section we will list (in table $2$) a representative seven-tuple $\bm$ for each
solution. (Recall that $\bm$ determines the overall conjugacy class of the
triple $M_1,M_2,M_3$.)
Thus from this data  one can for example easily compute the full permutation
representation of each solution curve as a cover of the three-punctured sphere.
(It is too cumbersome to write them all directly.)  
Rather than write the numbers $m_i=\tr(M_i)$ and $m_{ij}=\tr(M_iM_j)$ it is
simpler to write the rational numbers $\theta_i, \si_{ij}$, where
$$\tr(M_i)=2\cos(\pi \theta_i)\qquad\text{and}
\qquad \tr(M_iM_j)=2\cos(\pi \si_{ij})$$
with $0\leqslant\theta_i,\si_{ij}\leqslant 1$.

{\Small
\begin{table}
\begin{center}

\begin{tabular}[h]{|c|c|c||c|c|c| }
\hline
1& (1/2, 0, 1/3, 1/5) & (1/2, 1/3, 1/5) &\ \ 27& (2/5, 2/5, 2/3, 2/5) & (4/5, 3/5, 3/5) \\ \hline
2& (1/2, 0, 1/3, 2/5) & (1/2, 1/3, 2/5) &\ \ 28& (1/2, 1/2, 1/5, 3/5) & (1/2,1/2, 2/3) \\ \hline
3& (1/2, 0, 1/5, 2/5) & (1/2, 1/5, 2/5) &\ \ 29& (1/3, 1/3, 1/3, 4/5) & (2/3, 3/5, 3/5) \\ \hline
4& (1/3, 0, 1/5, 2/5) & (1/3, 1/5, 2/5) &\ \ 30& (1/3, 1/3, 1/3, 2/5) & (2/3, 1/5, 1/5) \\ \hline
5& (1/3, 0, 1/3, 1/5) & (1/3, 1/3, 1/5) &\ \ 31& (4/5, 4/5, 4/5, 4/5) & (1/3, 1/5, 0) \\ \hline
6& (1/3, 0, 1/3, 3/5) & (1/3, 1/3, 3/5) &\ \ 32& (3/5, 3/5, 3/5, 3/5) & (3/5, 3/5, 3/5) \\ \hline
7& (1/3, 0, 1/5, 1/5) & (1/3, 1/5, 1/5) &\ \ 33& (2/5, 1/2, 1/3, 4/5) & (4/5, 2/3, 2/3) \\ \hline
8& (1/3, 0, 2/5, 2/5) & (1/3, 2/5, 2/5) &\ \ 34& (1/5, 1/3, 1/5, 1/2) & (1/3, 2/5, 1/3) \\ \hline
9& (1/5, 0, 1/5, 1/5) & (1/5, 1/5, 1/5) &\ \ 35& (2/5, 1/3, 2/5, 1/2) & (1/5, 1/5, 4/5) \\ \hline
10& (3/5, 0, 3/5, 3/5) & (3/5, 3/5, 3/5) &\ \ 36& (1/3, 1/5, 1/3, 2/5) & (2/5, 2/5, 1/2) \\ \hline
11& (1/3, 1/5, 1/5, 2/3) & (1/2, 1/3, 1/2) &\ \ 37& (1/3, 1/3, 1/3, 1/5) & (1/5, 1/3, 1/2) \\ \hline
12& (1/3, 2/5, 2/5, 2/3) & (1/2, 1/3, 1/2) &\ \ 38& (1/3, 1/3, 1/3, 3/5) & (1/3, 1/3, 1/2) \\ \hline
13& (1/5, 2/5, 2/5, 4/5) & (1/2, 3/5, 1/2) &\ \ 39& (1/3, 4/5, 1/3, 4/5) & (2/3, 3/5, 0) \\ \hline
14& (2/5, 1/5, 1/5, 2/3) & (1/2, 1/3, 1/2) &\ \ 40& (3/5, 2/3, 3/5, 2/3) & (2/3, 1/5, 1/3) \\ \hline
15& (4/5, 2/5, 2/5, 2/3) & (1/2, 1/3, 1/2) &\ \ 41& (1/3, 1/3, 1/3, 1/3) & (1/3, 1/3, 3/5) \\ \hline
16& (1/5, 1/5, 1/5, 1/2) & (1/3, 2/5, 1/3) &\ \ 42& (1/3, 1/2, 1/3, 4/5) & (2/3, 4/5, 3/5) \\ \hline
17& (2/5, 2/5, 2/5, 1/2) & (1/3, 4/5, 1/3) &\ \ 43& (1/3, 1/2, 1/3, 2/5) & (2/3, 1/2, 1/3) \\ \hline
18& (1/5, 1/5, 1/5, 2/5) & (1/3, 1/5, 1/3) &\ \ 44& (1/2, 1/5, 1/2, 1/5) & (1/3, 2/3, 0) \\ \hline
19& (2/5, 2/5, 2/5, 1/5) & (1/3, 3/5, 1/3) &\ \ 45& (1/2, 2/5, 1/2, 2/5) & (1/5, 4/5, 0) \\ \hline
20& (2/5, 1/3, 1/5, 2/3) & (2/3, 1/3, 1/2) &\ \ 46& (1/3, 1/3, 1/3, 1/2) & (1/2, 1/3, 3/5) \\ \hline
21& (1/5, 1/5, 2/5, 2/5) & (1/3, 1/2,  1/3) &\ \ 47& (1/2, 1/2, 1/3, 1/5) & (1/2, 1/2, 3/5) \\ \hline
22& (1/5, 1/5, 2/5, 1/3) & (2/5, 1/3, 1/3) &\ \ 48& (1/2, 1/2, 1/3, 2/5) & (1/2, 4/5, 1/2) \\ \hline
23& (2/5, 2/5, 1/5, 2/3) & (4/5, 1/3, 1/3) &\ \ 49& (1/3, 1/2, 1/3, 1/2) & (2/3, 3/5, 3/5) \\ \hline
24& (1/2, 2/5, 1/5, 4/5) & (2/3, 1/3, 1/2) &\ \ 50& (1/2, 1/2, 1/2, 1/5) & (1/2, 2/5, 1/3) \\ \hline
25& (2/5, 2/5, 1/2, 4/5) & (3/5, 2/3, 1/2) &\ \ 51& (1/2, 1/2, 1/2, 2/5) & (1/2, 2/3, 1/5) \\ \hline
26& (1/5, 1/5, 1/5, 2/3) & (1/5, 2/5, 1/5) &\ \ 52& (1/2, 1/2, 1/2, 2/3) & (1/2, 1/5, 2/5) \\ \hline
\end{tabular}

%16, 17 changed by hand to fit into hitchin's k=3 (poncelet) dihedral family
%(checked puiseux at 0)

\vspace{0.2cm}
\caption{Representative seven-tuples $(\theta_1,\theta_2,\theta_3,\theta_4)$,
  $(\si_{12},\si_{23},\si_{13})$.}
\label{seventuple table}
\end{center}\end{table}
%\end{landscape}
}

\end{section}

\begin{section}{The generic icosahedral solution}\label{sn: generic}

Looking carefully at the `Walls' column of table $1$, 
the author was surprised to see there
is a zero, on row $33$. Namely there is a solution whose parameters
are generic in that they lie on {\em none} of the affine $F_4$ hyperplanes. 
(This clearly
implies they also lie on none of the reflecting hyperplanes of 
Okamoto's affine $D_4$ action.)

Being convinced there must be some mistake this solution was pursued
further and in this section we will present the explicit solution: 
an algebraic solution with $12$ branches and generic parameters; the largest
genus zero icosahedral solution. 
Note that this solution is generic in another sense too: it has the largest
value $2304$ of $n$ in table $1$, so a randomly chosen triple of generators of the
icosahedral group is most likely to lead to this solution.
(This number $2304=12\cdot 4!\cdot 2^3$ also shows that the group $\wmc$ is the
smallest possible extension of $\F_2$ yielding Theorem \ref{thm: exactly 52}: namely
the $\F_2$ orbits have size $12$ so the $2304$ points of $S$ correspond to 
$4!\cdot 2^3$ solutions of this type, each with $12$ branches, and so we deduce that the group 
$\wmc /\F_2\cong S_4\sdp \Si$ 
acts simply transitively on this set of solutions.)

Observe also, from the $A_5$ type, that the (projective) 
linear monodromy of the corresponding linear system 
has one generator in each
of the four non-trivial conjugacy class of $A_5$. 
In this 
sense we are at the opposite extreme to the Dubrovin--Mazzocco solutions,
which have all the $A_5$ conjugacy classes equal.

The solution was constructed as described in \cite{k2p}, using  Jimbo's
asymptotic formula (\cite{Jimbo82}, \cite{k2p} Theorem 4), noting 
that this solution has `good' representatives.
(The idea of using precise knowledge of the asymptotics to determine algebraic
solutions of $\text{P}_{\text{VI}}$ was also used in \cite{DubMaz00}.)
The solution curve thus arrived at is:
%\footnote{
%The first curve that was found was much longer---and a parameterisation looked
%beyond reach---but, by choosing a better
%representative, several months later the present curve was found.}

%\begin{table}[h]
%\begin{center}
%$$\text{\underline{Generic icosahedral solution}}$$
\begin{gather}
\text{\underline{Generic icosahedral solution}}\notag \\ 
%$(\theta_1,\theta_2,\theta_3,\theta_4)=(2/5, 1/2, 1/3, 4/5):$}\notag \\ 
 (15524784\,t^2-5373216\,t+1350000)\,y^{12}- 
 (128381760\,t^2-13366080\,t)\,y^{11}+ \notag \\
 (5425704\,t^3+496677744\,t^2-30539160\,t)\,y^{10}- \notag \\
 (14929920\,t^4+41364000\,t^3+866759680\,t^2-2928160\,t)\,y^9+ \notag \\
 (107546535\,t^4-508275750\,t^3+747613335\,t^2-1837080\,t)\,y^8- \notag \\
 (24385536\,t^5-285548724\,t^4-2437066824\,t^3+74927724\,t^2+944784\,t)\,y^7+ \notag \\
 (58212000\,t^5-2865570750\,t^4-4456260900\,t^3+17631810\,t^2)\,y^6- \notag \\
 (49787136\,t^6-904003584\,t^5-7215732804\,t^4-2130570936\,t^3-12872196\,t^2)\,y^5- \notag \\
 (413500320\,t^6+3724484160\,t^5+4839581265\,t^4+162430110\,t^3+3750705\,t^2)\,y^4+ \notag \\
 (3001304640\,t^6+74794560\,t^5+2710584000\,t^4-380946240\,t^3)\,y^3- \notag \\
 (940800000\,t^7+977540640\,t^6-726801696\,t^5+939255264\,t^4-72013536\,t^3)\,y^2+ \notag \\
 (1176000000\,t^7-1481095680\,t^6+765158400\,t^5)\,y- \notag \\
 (1920800000\,t^8-7212800000\,t^7+10522980864\,t^6-6913299456\,t^5+1728324864\,t^4) \notag
%(\theta_1,\theta_2,\theta_3,\theta_4)=(2/5, 1/2, 1/3, 4/5)\notag
\end{gather}
%$$(\theta_1,\theta_2,\theta_3,\theta_4)=(2/5, 1/2, 1/3, 4/5)$$
%\end{center}
%\end{table}

Implicit differentiation enables us to confirm that 
the function $y(t)$ defined by this polynomial is indeed a solution to
the Painlev\'e VI equation,
with $(\theta_1,\ldots\theta_4)=$
%$(2/5, 1/2, 1/3, 4/5)$. 
$(\frac{2}{5}, \frac{1}{2},\frac{1}{3},\frac{4}{5})$.
Since this represents 
a genus zero curve one can use a computer to find a rational
parameterisation. The author is very grateful to Mark van Hoeij for performing
this task for the above curve and finding that it may be parameterised simply, 
as in Theorem B in the introduction.

%$$y=-\frac{9 s (s^2+1) (3 s-4) (15 s^4-5 s^3+3 s^2-3 s+2)}
%{(2 s-1)^2(9 s^2+4)(9 s^2+3 s+10)},\qquad 
%t=\frac{27 s^5 (s^2+1)^2 (3 s-4)^3}{4(2 s-1)^3(9 s^2+4)^2}.$$

\end{section}

\begin{section}{More examples}\label{sn: moreex}

By now it is not too much 
extra trouble to produce other solutions (since the
procedure was sufficiently systematised in order to
find a simple version of the generic solution). Here are the remaining 
good genus zero icosahedral solutions:

\begin{gather}
\text{Solution $20$, genus zero, $5$ branches,
$(\theta_1,\theta_2,\theta_3,\theta_4)=(2/5, 1/3, 1/5, 2/3)$:}\notag\\
%has rep with all sigmas = 1/2
y=
{\frac {2 \left ({s}^{2}+s+7\right )\left (5\,s-2\right )}{s\left (s+
5\right )\left (4\,{s}^{2}-5\,s+10\right )}},\qquad 
t={\frac {27 \left (5\,s-2\right )^{2}}{\left (s+5\right )\left (4\,{s}
^{2}-5\,s+10\right )^{2}}}\notag
\end{gather}

\begin{gather}
\text{Solution $24$, genus zero, $8$ branches, 
$(\theta_1,\theta_2,\theta_3,\theta_4)=(1/2, 2/5, 1/5, 4/5)$:}\notag \\
%has rep with sigmas = 2/3, 1/3, 1/2
y={\frac {s\left (s+4\right )\left (3\,{s}^{4}-2\,{s}^{3}-2\,{s}^{2
}+8\,s+8\right )}{8 \left (s-1\right )\left ({s}^{2}+4\right )\left (s+1
\right )^{2}}},\qquad
t={\frac {{s}^{5}\left (s+4\right )^{3}}
{4 \left (s-1\right )\left ({s}^{2}+4\right )^{2}\left (s+1\right )^{3}}}
\notag
\end{gather}

\begin{gather}
\text{Solution $25$, genus zero, $8$ branches,
$(\theta_1,\theta_2,\theta_3,\theta_4)=(2/5, 2/5, 1/2, 4/5)$:}\notag \\
%has rep with sigmas = 1/3, 1/2, 1/2
y={\frac {{s}^{2}\left (5\,{s}^{3}+2\,{s}^{2}-4\,s-8\right )\left (
s+4\right )^{2}}
{4 \left (s+1\right )^{2}\left ({s}^{2}+4\right )\left (
s-1\right )\left ({s}^{2}+3\,s+6\right )}},\qquad
t={\frac {{s}^{5}\left (s+4\right )^{3}}
{4 \left (s-1\right )\left ({s}^{2}+4\right )^{2}\left (s+1\right )^{3}}}
\notag
\end{gather}

\begin{gather}
\text{Solution $28$, genus zero, $10$ branches, 
$(\theta_1,\theta_2,\theta_3,\theta_4)=(1/2, 1/2, 1/5,  3/5)$:} \notag \\
%has rep with all sigmas = 2/3
y={\frac { \left( {s}^{5}+5\,{s}^{4}-20\,{s}^{3}+75\,s+75 \right) 
 \left( {s}^{2}-5 \right)  \left( {s}^{2}+5 \right) }{ \left( s+1
 \right) ^{2} \left( {s}^{2}-4\,s+5 \right)  \left( s+5 \right) 
 \left( {s}^{4}+6\,{s}^{2}-75 \right) }},\quad
t=
{\frac { 2 \left( {s}^{2}+5 \right) ^{3} \left( {s}^{2}-5 \right) ^{2
}}{ \left( s+5 \right) ^{3} \left( {s}^{2}-4\,s+5 \right) ^{2} \left( 
s+1 \right) ^{3}}}
\notag
\end{gather}

The next solution we can find using Jimbo's formula 
is the generic solution (number $33$) already
displayed.
Beyond that we pass onto the higher genus solutions. 
In principle 
we can still find these, although 
eventually one will have trouble computing all the symmetric functions of the
Puiseux series of the solutions on the branches.

For example solutions $34$ and $35$ both 
become single valued on the elliptic curve:
\beq\label{eq: ell curve}
u^2=(3s+5)(8s^2-5s+5)
\eeq
and, as functions on this curve, the solutions are 
given explicitly as:
\begin{gather}
\text{Solution $34$, genus one, $12$ branches, 
$(\theta_1,\theta_2,\theta_3,\theta_4)=(1/5, 1/3, 1/5, 1/2)$:}\notag \\
y=\frac{1}{2}+
{\frac { \left( 3\,s+5 \right)  \left( 8\,{s}^{4}-10\,{s}^{3}+12
\,{s}^{2}-13\,s+11 \right) }{ 2\,\left( 2\,{s}^{3}-15\,s+5 \right) u}}
\notag \\
t=t_{34}=
\frac{1}{2}-{\frac { \left( 8\,{s}^{6}+20\,{s}^{3}-15\,{s}^{2}+66\,s-15
 \right)}{ 2\,\left( 8\,{s}^{2}-5\,s+5 \right) \, u}}\notag
\end{gather}

\begin{gather}
\text{Solution $35$, genus one, $12$ branches, 
$(\theta_1,\theta_2,\theta_3,\theta_4)=(2/5, 1/3, 2/5, 1/2)$:}\notag \\
y=\frac{1}{2}+
{\frac { \left( 3\,s+5 \right)  \left( 16\,{s}^{5}-8\,{s}^{4}+18
\,{s}^{3}-8\,{s}^{2}+115\,s+3 \right) }{2\,\left( 26\,{s}^{3}+60\,{s}^{
2}+15\,s+35 \right)\, u }},\quad t=t_{34}
\notag
\end{gather}

Next, solution $36$ is given by the functions below on the curve 
$u^2=3(5s+1)(8s^2-9s+3)$:
\begin{gather}
\text{Solution $36$, genus one, $12$ branches, 
$(\theta_1,\theta_2,\theta_3,\theta_4)=(1/3, 1/5, 1/3, 2/5)$:}\notag \\
y=\frac{1}{2}+
{\frac {140\,{s}^{6}+1029\,{s}^{5}-1023\,{s}^{4}+360\,{s}^{3}-
288\,{s}^{2}+27\,s+27}{18\,u \left( s+1 \right)  \left( 7\,{s}^{3}-3\,{s}^
{2}-s+1 \right) }}\notag\\
t=\frac{1}{2}+
{\frac {40\,{s}^{6}+540\,{s}^{5}-765\,{s}^{4}+540\,{s}^{3}-270\,{
s}^{2}+27}{6\,u \left( 8\,{s}^{2}-9\,s+3 \right)  \left( s+1 \right) ^{2}
}}
\notag
\end{gather}
%where $u^2=3(5s+1)(8s^2-9s+3)$.

The next two solutions are related to the Valentiner group and will appear in
section \ref{sn: valentiner} and the outstanding good 
solutions with fewer than $20$
branches are:

\begin{gather}
\text{Solution $39$, genus one, $15$ branches, 
$(\theta_1,\theta_2,\theta_3,\theta_4)=(1/3, 4/5, 1/3, 4/5)$:}\notag \\
y=\frac{1}{2}+
{\frac {14\,{s}^{5}+61\,{s}^{4}-66\,{s}^{3}-660\,{s}^{2}-900\,s-
225}{ 6\,\left( s+1 \right)  \left( {s}^{2}-5 \right) u}}
\notag\\
t=\frac{1}{2}-
{\frac { \left( 2\,{s}^{7}+10\,{s}^{6}-90\,{s}^{4}-135\,{s}^{3}
+297\,{s}^{2}+945\,s+675 \right) u}{18\, \left( 4\,{s}^{2}+15\,s+15
 \right) ^{2} \left( {s}^{2}-5 \right) }}
\notag
\end{gather}
where
${u}^{2}=3\, \left( s+5 \right)  \left( 4\,{s}^{2}+15\,s+15 \right).$

\begin{gather}
\text{Solution $40$, genus one, $15$ branches, 
$(\theta_1,\theta_2,\theta_3,\theta_4)=(3/5, 2/3, 3/5, 2/3)$:}\notag \\
y=\frac{1}{2}
-{\frac {2\,{s}^{9}+20\,{s}^{8}+53\,{s}^{7}-89\,{s}^{6}-605\,{s}^
{5}-851\,{s}^{4}-1389\,{s}^{3}-5775\,{s}^{2}-10125\,s-5625}{2\, \left( {s
}^{2}-5 \right)  \left( {s}^{2}-6\,s-15 \right)  \left( {s}^{2}+4\,s+5
 \right) u}}
\notag
\end{gather}
where $t,u,s$ are as for solution $39$ above.

Now we will fill in the gaps and explain how one may
find the outstanding solutions 
for which 
Jimbo's formula cannot  be applied directly, namely no.s 
$22,23,27,29,30$.
Upon inspection one finds that these solutions always have a regular
branch at zero (namely there is a cycle of lenth one in the 
permutation of the branches of the solution curve at zero).
Thus we need to find the leading term in the Taylor/Laurent expansion of the
solutions on the regular branches (the leading terms of 
the Puiseux expansions on the other branches still 
being given by Jimbo's formula). 
First we observe that each of these five solutions has a representative for
which
$$\theta_1+\theta_2=\si$$
on one branch at zero (the regular branch),
where $2\cos(\pi\si)=m_{12}=\tr(M_1M_2),\, 0<\re(\si)< 1$.
Then the leading term is given by the following result.
\begin{lem}
If $\theta_1+\theta_2=\si$ on a branch of a solution to Painlev\'e VI
with finite linear monodromy group, then the leading term of the Laurent
expansion at zero of the solution is
$$y(t)= \frac{\theta_1}{\theta_1+\theta_2} t +O(t^2)$$
\end{lem}
\sketch
We proceed as in Jimbo's article \cite{Jimbo82};
As $t\to 0$ the system \eqref{eq: lin syst}
degenerates into two hypergeometric systems (\cite{Jimbo82} (2.13, 2.14))
and the fundamental solutions and monodromy data can be related explicitly.
Solving the Riemann--Hilbert problems for the two hypergeometric systems
gives the asymptotics
for the isomonodromic family of systems 
\eqref{eq: lin syst} (see \cite{Jimbo82} (2.15)), 
and therefore also for the $\text{P}_{\text{VI}}$ solution.

In the case we are considering the condition $\theta_1+\theta_2=\si$ 
forces one of the hypergeometric systems (\cite{Jimbo82} (2.14)) to be
reducible, and in fact abelian---since the monodromy group is finite.
This makes the corresponding Riemann--Hilbert problem very easy to solve 
and yields the stated formula for the leading term.
\esketch

\begin{rmk}
This can almost be guessed directly:
substituting $y=a_1t+a_2t^2$ into $\text{P}_{\text{VI}}$ gives the leading term
$$\frac{((\theta_2+\theta_1)a_1-\theta_1)((\theta_2-\theta_1)a_1+\theta_1)}
{2a_1(1-a_1)}   t^{-1} 
$$
at zero.
Thus if $\theta_2=\theta_1$ then
the value of $a_1$ in the lemma is forced, and moreover 
our five examples all have representatives with
$\theta_2=\theta_1$.
\end{rmk}

The explicit formulae for these five solutions thus obtained are as
follows. 

\begin{gather}
\text{Solution $22$, genus zero, $6$ branches, 
$(\theta_1,\theta_2,\theta_3,\theta_4)=(1/5, 1/5, 2/5, 1/3)$:}\notag \\
y={\frac { -54\,s\left( s-7 \right)}{ \left({s}^{4} -20\,{s}^{2}-35 \right)  
\left( s+1 \right)  \left( s-4 \right) }},\qquad
t=
{\frac {432\,s}
{ \left( s+5 \right)  \left( s-4 \right) ^{2} \left( s+1 \right) ^{3}}}
\notag
\end{gather}

\begin{gather}
\text{Solution $23$, genus zero, $6$ branches, 
$(\theta_1,\theta_2,\theta_3,\theta_4)=(2/5, 2/5, 1/5, 2/3)$:}\notag \\
y=
{\frac { 18\,s\,\left( s-3 \right)}
{ \left( s-4 \right)  \left( s+1 \right)  \left( {s}^{2}+5 \right) }},
\qquad
t=
{\frac {432\,s}
{ \left( s+5 \right)  \left( s-4 \right) ^{2} \left( s+1 \right) ^{3}}}
\notag
\end{gather}

\begin{gather}
\text{Solution $27$, genus one, $9$ branches, 
$(\theta_1,\theta_2,\theta_3,\theta_4)=(2/5, 2/5, 2/3, 2/5)$:}\notag \\
y=\frac{1}{2}+{\frac { \left( 350\,{s}^{3}+63\,{s}^{2}-6\,s-2 \right)}
{30\,\left( 2\,s+1 \right) {s} u}}\quad
t=\frac{1}{2}+{\frac { \left( 25\,{s}^{4}+170\,{s}^{3}+42\,{s}^
{2}+8\,s-2 \right) u}{54\, \left( 5\,s+4 \right) ^{2}{s}^{3}}}\notag
\end{gather}

\noindent
where $u$ and $s$ satisfy 
$u^2=s \left( 8\,s+1 \right)  \left( 5\,s+4 \right).$

\begin{gather}
\text{Solution $29$, genus zero, $10$ branches, 
$(\theta_1,\theta_2,\theta_3,\theta_4)=(1/3, 1/3, 1/3, 4/5)$:}\notag \\
y=
{\frac { \left( s+2 \right)  \left( {s}^{2}+1 \right)  
\left(   2\,{s}^{2}+3\,s+3 \right) {s}^{2}}
{ 2\,\left( {s}^{2}+s+1 \right)  \left( 3\,{s}^{2}+3\,s+2
 \right) }},
\qquad
t=t_{29}=
{\frac { \left( s+2 \right)  \left( 2\,{s}^{2}+3\,s+3 \right) ^{2}{s}^{5}}{ \left( 2\,s+1 \right)  \left( 3\,{s}^{2}+3\,s+2 \right) ^{2}}}
\notag
\end{gather}

\begin{gather}
\text{Solution $30$, genus zero, $10$ branches, 
$(\theta_1,\theta_2,\theta_3,\theta_4)=(1/3,1/3,1/3,2/5)$:}\notag \\
y=
{\frac { \left( s+2 \right)  \left( 2\,{s}^{2}+3\,s+3 \right)  \left( 7\,{s}^{2}+10\,s+7 \right) {s}^{4}}{ \left( 3\,{s}^{2}+3\,s+2 \right)  \left( 4\,{s}^
{6}+12\,{s}^{5}+15\,{s}^{4}+10\,{s}^{3}+15\,{s}^{2}+12\,s+4 \right) }},
\quad
t=t_{29}
\notag
\end{gather}

%\newpage

\begin{rmk} (Added May 2005.)\     \label{rmk: qt}
The remaining solutions (except the Valentiner solutions, 
which will appear in 
section \ref{sn: valentiner} below, and solutions 42 and 43) 
may be obtained from known solutions using the quadratic
transformations defined in 1991 by Kitaev \cite{Kitaev-quad-p6}.
The basic idea is as follows. 
Given an icosahedral Fuchsian system $A$ with $A_5$ type
$a^2\xi\eta$ for some $\xi,\eta\in\{a,b,c,d\}$
(i.e. with two local monodromies, say at $0$ and $\infty$, of order two in
$\PSL_2(\IC)$)
we can pull back along the map $w\mapsto z=w^2$,
and remove the resulting apparent singularities, 
to get a Fuchsian 
%system with two apparent singularities at $0$, $\infty$
%and four non-apparent singularities at $\pm 1,\pm\sqrt{t}$.
%Removing the apparent singularities (using Schlesinger transformations)
%yields a 
system $B$ with $A_5$ type $\xi^2\eta^2$.
%, which may be put in the form
%\eqref{eq: lin syst} by a coordinate transformation.
Isomonodromic deformations of $A$ correspond to isomonodromic deformations of
$B$, and one can obtain formulae relating the corresponding \PVI\ solutions.
In practice the formulae are much simpler at different (Okamoto
equivalent) values of the parameters (see \cite{RGT-quad} (2.7) and also
the recent article \cite{TOS-folding}).
In the cases at hand, this procedure gives an algebraic relation 
with a solution having half the number of branches;
Examining Table 1 we see solution 31 $\Rightarrow$ solution 44 and in turn
solution 44 $\Rightarrow$ solution 50.
Similarly 
%32 $\Rightarrow$ 45 $\Rightarrow$ 51,
%39 $\Rightarrow$ 47,
%40 $\Rightarrow$ 48,
%and 41 $\Rightarrow$ 49 $\Rightarrow$ 52.
$$
32  \Rightarrow  45  \Rightarrow  51,\qquad
39  \Rightarrow  47,\qquad
40  \Rightarrow  48,\qquad
41  \Rightarrow  49  \Rightarrow  52.$$
(Some work is still required to obtain efficient parameterisations of the
solutions obtained in this way.)
Fortunately the only possible simplification to the rest of this article is 
21 $\Rightarrow$ 28, which was no trouble anyway.

In any case this gave us motivation to construct the remaining 
solutions 42 and 43 
using our original method, essentially completing the construction 
of all icosahedral solutions:

\begin{gather}
\text{Solution $42$, genus one, $20$ branches, 
$(\theta_1,\theta_2,\theta_3,\theta_4)=(1/3, 1/2, 1/3, 4/5)$:}\notag \\
y=\frac{1}{2}+
{\frac { \left( 8\,{s}^{6}-28\,{s}^{5}+85\,{s}^{4}-196\,{s}^{3}+
214\,{s}^{2}-196\,s+41 \right)  \left( s+3 \right) }
{6\left( {s}^{2}+1 \right)  \left( 3\,{s}^{2}-4\,s+5 \right)u }}
\notag\\
t=\frac{1}{2}-
{\frac { \left( s+3 \right) \text{P}
%  \left( 8\,{s}^{10}+100\,{s}^{7}-135
%\,{s}^{6}+834\,{s}^{5}-1205\,{s}^{4}+2280\,{s}^{3}-1365\,{s}^{2}+890\,
%s+321 \right) 
}
{2 \left( {s}^{2}+1 \right) ^{2}{u}^{3}}}
\notag
\end{gather}
where
${u}^{2}=3\left( s+3 \right)\left( 8\,{s}^{2}-13\,s+17 \right),$ and
$$\text{P} = 8\,{s}^{10}+100\,{s}^{7}-135
\,{s}^{6}+834\,{s}^{5}-1205\,{s}^{4}+2280\,{s}^{3}-1365\,{s}^{2}+890\,
s+321.$$

\begin{gather}
\text{Solution $43$, genus one, $20$ branches, 
$(\theta_1,\theta_2,\theta_3,\theta_4)=(1/3, 1/2, 1/3, 2/5)$:}\notag \\
y=\frac{1}{2}+
{\frac { \left( s+3 \right)\text{Q}
%\left( 28\,{s}^{9}-235\,{s}^{8}+556\,{s}^{7}-1334\,{s}^
%{6}+2174\,{s}^{5}-3854\,{s}^{4}+4360\,{s}^{3}-4738\,{s}^{2}+2362\,s-
%1047 \right)  
}
{ 18 \left( {s}^{2}+1 \right)  
\left( {s}^{6}-7\,{s}^{4}+42\,{s}^{3}-45\,{s}^{2}+34\,s+7 \right) u}}
\notag
\end{gather}
where
$$\text{Q}=28\,{s}^{9}-235\,{s}^{8}+556\,{s}^{7}-1334\,{s}^
{6}+2174\,{s}^{5}-3854\,{s}^{4}+4360\,{s}^{3}-4738\,{s}^{2}+2362\,s-
1047$$
and $t,u,s$ are as for solution 42.
\end{rmk}

%\begin{rmk}
%The elliptic curves appearing here (in solutions $27,34,36$) all have rank
%zero, and so have a finite number of rational points (making up groups
%$\IZ/2\times \IZ/6, \IZ/10, \IZ/12$ resp.).
%Thus there appear to be
%relatively few systems of differential equations with $\IQ$
%coefficients having the corresponding monodromy representations.
%\end{rmk}

\end{section}

\begin{section}{The Valentiner Solutions}\label{sn: valentiner}

The Valentiner reflection group is the subgroup of $\GL_3(\IC)$
generated by the complex reflections (cf. e.g. \cite{Shep-Todd}):
$$
r_1= \left( \begin {array}{ccc} 0&-{\omega}^{2}&0\\\noalign{\medskip}
-\omega&0&0\\\noalign{\medskip}0&0&1\end {array} \right) 
\quad
r_2=-\frac{1}{2}
 \left( \begin {array}{ccc} -1 &\omega\tau&{\frac {\omega^2}{
\tau}}\\\noalign{\medskip}{\frac {\tau}{\omega}}&
{\tau}^{-1}&\omega\\\noalign{\medskip}{\frac {
\omega}{\tau}}&{\omega}^{2}&-\tau\end {array} \right) 
\quad
r_3= \left( \begin {array}{ccc} -1&0&0\\\noalign{\medskip}0&1&0
\\\noalign{\medskip}0&0&1\end {array} \right) 
$$
where $\omega=\exp(2\pi i/3)$, $\tau=(1+\sqrt{5})/2$.
It has order $2160$ and
the corresponding projective group in $\PGL_3(\IC)$ is isomorphic to $A_6$ the
alternating group on six letters.

We wish to apply the procedure of \cite{k2p} section 2 
to this triple of generating
reflections to obtain a triple of elements of $\SL_2(\IC)$.
By definition the corresponding $\SL_2(\IC)$ triple $(M_1,M_2,M_3)$
has invariants
\begin{gather}
\tr(M_1)=\frac{t_1}{n_1}+\frac{n_1}{t_1},\qquad
\tr(M_2)=\frac{t_2}{n_1}+\frac{n_1}{t_2},\qquad
\tr(M_3)=\frac{t_3}{n_1}+\frac{n_1}{t_3}, \notag \\
\tr(M_1M_2)= \frac{t_{12}}{t_1t_2},\qquad
\tr(M_2M_3)= \frac{t_{23}}{t_2t_3},\qquad
\tr(M_1M_3)= \frac{t_{13}}{t_1t_3},\notag \\
\tr(M_4)=\tr(M_3M_2M_1)=\frac{n_2}{n_3}+\frac{n_3}{n_2} \notag
\end{gather}
where 
$t_{jk}=\tr(r_jr_k)-1$, $t_j$ is a choice of square root of $\det(r_j)$,
and the $n_j$ are chosen square roots of the eigenvalues of 
the product $r_3r_2r_1$ (which we are thus choosing an order of too).
Here each of the reflections $r_j$ is of order two
so we take may take the invariant $t_j=i$ for each $j$.
Next, the product $r_3r_2r_1$ has eigenvalues 
$\{
\exp(2\pi i\frac{5}{30}),
\exp(2\pi i\frac{11}{30}),
\exp(2\pi i\frac{29}{30}) \}$,
so we may take 
$$
n_1=\exp(5\pi i/{30}),\quad
n_2=\exp(11\pi i/{30}),\quad
n_3=\exp(29\pi i/{30}).$$
Also we compute:
$$\tr(r_1r_2)=0,\qquad
\tr(r_2r_3)=0,\qquad
\tr(r_1r_3)=1.$$
Then the corresponding $\SL_2(\IC)$ invariants are:
$$m_1=m_2=m_3=1,
\quad m_{4}=2\cos(3\pi/5),\quad 
m_{12}=m_{23}=1,\quad m_{13}=0.$$
Thus the $\theta$ parameters are $(1/3,1/3,1/3,3/5)$,
since $2\cos(\pi/3)=1$, and 
one finds then that the corresponding $\SL_2(\IC)$ triple
generates the binary icosahedral group 
and corresponds to row $38$ of tables $1,2$.

In particular if we are able to find the corresponding \PVI\ solution then (as
was done in \cite{k2p} for the Klein group) we
can explicitly construct an isomonodromic family of rank three Fuchsian
equations (with four poles on $\IP^1$) having monodromy group equal to the
Valentiner reflection group, generated by reflections.

Rather than repeat the details (which are exactly as in \cite{k2p}) we just
give the \PVI\ solution:

\begin{gather}
\text{Solution $38$, genus one, $15$ branches, 
$(\theta_1,\theta_2,\theta_3,\theta_4)=(1/3, 1/3, 1/3, 3/5)$:}\notag \\
y=\frac{1}{2}
-{\frac {250\,{s}^{6}+500\,{s}^{5}+518\,{s}^{4}+261\,{s}^{3}+76\,
{s}^{2}+13\,s+2}{2\, \left( s+2 \right)  \left( 5\,s+1 \right)  \left( 5
\,{s}^{3}+6\,{s}^{2}+3\,s+1 \right) u}}\notag\\
t=\frac{1}{2}
-{\frac {3\,(500\,{s}^{7}+925\,{s}^{6}+1164\,{s}^{5}+830\,{s}^{4}+340
\,{s}^{3}+105\,{s}^{2}+20\,s+4)}{ 2\,\left( s+2 \right) ^{2} \left( 5\,s+1
 \right) {u}^{3}}}
\notag
\end{gather}
where $(u,s)$ lies on the elliptic curve 
$${u}^{2}=\left( 4\,{s}^{2}+s+1 \right)  \left( 5\,s+1 \right). $$

However, unlike in the case of the Klein reflection group, the 
Valentiner group
has three inequivalent triples of generating reflections (above we used the
standard generating triple whose product has eigenvalues involving the
exponents of the group).
This is similar to the case of the icosahedral reflection group studied in 
\cite{DubMaz00}, which also has three inequivalent triples of generating
reflections (leading to the icosahedral solutions on rows $31,32,41$) 
although
now all three solutions are elliptic and the largest has $24$ branches.

The second generating triple gives the sibling solution of that above and
arises by replacing $r_1$ above by 
$$r_1=\left(\begin{matrix}
0&0&-{\omega}\\
0&1&0\\
-\omega^2&0&0\end{matrix} \right). $$
%which is the unitary reflection negating $e_1+\omega e_3$ so is in
%G_2160 
Then the product $r_3r_2r_1$ has eigenvalues 
$
\exp(2\pi i\frac{5}{30}),$
$\exp(2\pi i\frac{17}{30}),$
$\exp(2\pi i\frac{23}{30})$, and similarly to
 above one finds the correpsonding $\SL_2(\IC)$
triple corresponds to row $37$ of tables $1$ and $2$.
The corresponding \PVI\ solution is:

\begin{gather}
\text{Solution $37$, genus one, $15$ branches, 
$(\theta_1,\theta_2,\theta_3,\theta_4)=(1/3, 1/3, 1/3, 1/5)$:}\notag \\
y=\frac{1}{2}
-{\frac {1000\,{s}^{8}+2425\,{s}^{7}+4171\,{s}^{6}+3805\,{s}^{5}+
1999\,{s}^{4}+874\,{s}^{3}+244\,{s}^{2}+58\,s+4}{4\, \left( s+2 \right) 
 \left( 25\,{s}^{6}+135\,{s}^{5}+111\,{s}^{4}+91\,{s}^{3}+36\,{s}^{2}+
6\,s+1 \right) u}}
\notag
\end{gather}
with $t,u,s$ as for solution $38$ above.

Finally the third generating triple of the Valentiner reflection group 
arises by replacing $r_1$ above by the reflection
$$
r_1=
\frac{1}{2}
\left(
\begin{matrix}
1-\tau & \tau & 1 \\
\tau & 1 & 1-\tau \\
1 & 1-\tau & \tau 
\end{matrix}
\right).
$$
Then the product $r_3r_2r_1$ has eigenvalues 
$\exp(2\pi i\frac{2}{12}),$ 
$\exp(2\pi i\frac{5}{12}),$ 
$\exp(2\pi i\frac{11}{12})$, and similarly to
 above one now finds the correpsonding $\SL_2(\IC)$
triple corresponds to row $46$ of tables $1$ and $2$.
The corresponding \PVI\ solution has $24$ branches so is larger than any
previously constructed solution (and currently there are no known elliptic
solutions of higher degree, regardless of whether they have been 
explicitly constructed). 
For this solution 
the previous method of constructing the solution polynomial, involving
computing the symmetric functions of the Puiseux expansion at $0$ 
of the solution branches, no longer works. 
(For example in the worst case one faces a sum of 
$\left(\!\begin{smallmatrix} 24\\12 \end{smallmatrix}\!\right)$ 
terms, each 
of which is a $12$-fold product of Puiseux expansions with many terms.)
Instead one can obtain some coefficients in this way and then use 
the expected Okamoto symmetries of the solution to determine the outstanding
coefficients, by solving some sparse overdetermined linear equations.
(One then checks by implicit differentiation that the resulting polynomial
indeed defines a \PVI\ solution.) 
Using Mark van Hoeij's algorithms (in the Maple algebraic curves
package) we then obtain the following parameterisation:

\begin{gather}
\text{Solution $46$, genus one, $24$ branches,
$(\theta_1,\theta_2,\theta_3,\theta_4)=(1/3, 1/3, 1/3, 1/2)$} \notag\\
y=\frac{1}{2}
-\frac{\text{P}}
{{2\, \left( 3\,{s}^{2}-2\,s+2 \right)  \text{R} \,u}},
\quad
t=
\frac{1}{2}
+{\frac { \left( {s}^{2}+4\,s-2 \right)  \text{Q}}{ 2\,\left( s+2 \right) 
 \left( 3\,{s}^{2}-2\,s+2 \right) ^{2}{u}^{3}}}
\notag
\end{gather}
where
{\Small
$$\text{P}=16\,{s}^{11}+72\,{s}^{10}+50\,{s}^{9}-242\,{s}^{8}-
3143\,{s}^{7}+6562\,{s}^{6}-8312\,{s}^{5}+9760\,{s}^{4}-9836\,{s}^{3}+
6216\,{s}^{2}-2288\,s+416,$$
$$\text{Q}=8\,{s}^{10}+16
\,{s}^{9}+24\,{s}^{8}-84\,{s}^{7}+429\,{s}^{6}-312\,{s}^{5}+258\,{s}^{
4}-288\,{s}^{3}+288\,{s}^{2}-128\,s+32,$$
$$\text{R}
=26\,{s}^{6}+18\,{s}^{5}-75\,{s}^{4}+50\,{s}^{3}+270\,{s}^{2}-312\,s+
104,$$
}
and where $(u,s)$ lies on the elliptic curve 
$$u^2=\left( 8\,{s}^{2}-7\,s+2 \right)  \left( s+2 \right).$$

\end{section}

\renewcommand{\baselinestretch}{1}              %
\normalsize
\bibliographystyle{amsplain}    \label{biby}
\bibliography{../thesis/syr}    
\end{document}